\documentclass[12pt]{article}

\usepackage[english]{babel}

\usepackage{amsfonts}
\usepackage{amssymb}
\usepackage{amsmath}
\usepackage{amsthm}
\usepackage{epsf}
\usepackage{float}

\usepackage{graphicx}

\usepackage{cite}

\usepackage{fullpage}

\newtheorem{theorem}{Theorem}[section]
\newtheorem{lemma}[theorem]{Lemma}
\newtheorem{condition}[theorem]{Condition}
\newtheorem{prp}[theorem]{Proposition}
\newtheorem{requirement}{Requirement}

\theoremstyle{definition}
\newtheorem{definition}[theorem]{Definition}
\newtheorem{remark}[theorem]{Remark}

\numberwithin{equation}{section} \numberwithin{figure}{section}

\newcommand{\const}{{\rm const}}

\renewcommand{\phi}{{\varphi}}

\newcommand{\cH}{{\mathcal H}}

\newcommand{\bF}{{\mathbf F}}

\newcommand{\bbR}{{\mathbb R}}
\newcommand{\bbN}{{\mathbb N}}

\newcommand{\bbZ}{{\mathbb Z}}

\newcommand{\Jnk}{J_{n, k}}
\newcommand{\Jnkm}{J_{n, k}^{\rm main}}

\newcommand{\amin}{a^{\rm min}}
\newcommand{\amax}{a^{\rm max}}

\let\phi=\varphi

\newcommand{\ft}{\tilde{f}}

\newcommand{\At}{\tilde{A}}
\newcommand{\Ft}{\tilde{F}}

\newcommand{\ep}{\varepsilon}
\newcommand{\al}{\alpha}
\newcommand{\be}{\beta}
\newcommand{\sref}[1]{(\ref{#1})}
\newcommand{\lr}[1]{\!\left(\!#1\!\right)}
\newcommand{\lra}[1]{\left|#1\right|}
\newcommand{\str}[1]{\textrm{$#1$}}
\newcommand{\dx}{dx}

\newcommand{\erm}[1]{\textrm{$#1$}}

\title{Spatially discrete reaction-diffusion equations with discontinuous hysteresis}
\author{Pavel  Gurevich\footnote{Free University of Berlin, RUDN University, email: gurevich@math.fu-berlin.de},
        Sergey  Tikhomirov\footnote{Saint-Petersburg State Univeristy; email: s.tikhomirov@spbu.ru}}

\begin{document}

\date{}

\maketitle

\vspace{-25pt}

\begin{abstract}
We address the question: Why may reaction-diffusion equations with hysteretic nonlinearities become ill-posed and how to amend this? To do so, we discretize the spatial variable and obtain a lattice dynamical system with a hysteretic nonlinearity. We analyze a new mechanism that leads to appearance of a spatio-temporal pattern called {\it rattling}: the solution exhibits a propagation phenomenon different from the classical traveling wave, while the hysteretic nonlinearity, loosely speaking, takes a different value at every second spatial point, independently of the grid size. Such a dynamics indicates how one should redefine hysteresis to make the continuous problem well-posed and how the solution will then behave. In the present paper, we develop main tools for the analysis of the spatially discrete model and apply them to a prototype case. In particular, we prove that the propagation velocity is of order $a t^{-1/2}$ as $t\to\infty$ and explicitly find the rate $a$.
\end{abstract}

\tableofcontents

\section{Introduction}

\subsection{Background}\label{subsecBackground}

Hysteresis, or, more generally, bistability, refers to a class of
nonlinear phenomena which are observed in numerous real-world
systems. It arises in description of ferromagnetic materials,
shape-memory alloys, elasto-plastic bodies, as well as many
biological, economical, and social models,
see~\cite{KrasnBook,VisintinBook,BrokateBook,MayergoyzBook,MielkeHysteresis,
KrejciBook}. The primary goal of the present paper is to analyze a new
mechanism (which we call {\em rattling}) for pattern formation in
spatially discrete systems of reaction-diffusion equations (lattice
dynamical systems) with hysteresis. The phenomenon occurs in any
space dimension, including dimension one, and persists even for
scalar equations. As it is explained below,  our results are
relevant not only for lattice dynamical systems, but also for
continuous systems with hysteresis. On the other hand, they link
pattern formation mechanisms in hysteretic and bistable slow-fast
systems.

Let us begin with the prototype spatially continuous problem
\begin{equation}\label{eqPrototype}
\left\{
\begin{aligned}
& v_\tau  =   v_{xx}+\cH(v),\quad x\in(-1,1),\ \tau>0,\\
& v(x,0)  =\varphi(x),\quad x\in(-1,1),
\end{aligned}\right.
\end{equation}
supplemented with, e.g., Neumann boundary conditions. Here
$\cH(\cdot)$ is the simplest {\em hysteresis operator}, namely, the
{\em non-ideal relay} or {\em bistable switch}, see
Fig.~\ref{figRelay}.a and the (slightly modified) rigorous
definition in Section~\ref{secSetting}.

\noindent
\begin{minipage}{\linewidth}
      \centering
\begin{minipage}{0.67\linewidth}
    \begin{figure}[H]
    {\ \hfill\epsfxsize112mm\epsfbox{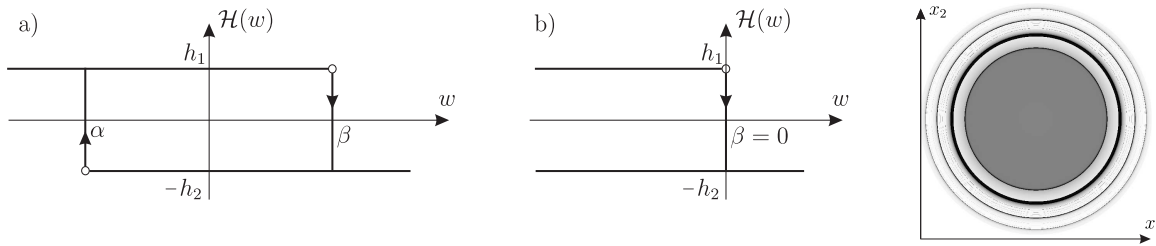}\hfill\ }  \caption{Hysteresis a)
     with thresholds $\alpha<\beta$. b) with
    thresholds $\alpha=-\infty$ and $\beta=0$.} \label{figRelay}
    \end{figure}
\end{minipage}
      \hspace{0.01\linewidth}
\begin{minipage}{0.28\linewidth}
    \begin{figure}[H]
    {\ \hfill\epsfxsize33mm\epsfbox{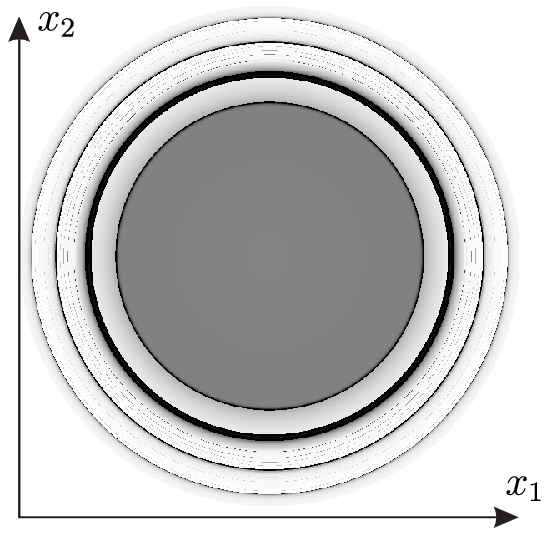}\hfill\ }
    \caption{Bacteria density at end of experiment.} \label{figRings}
    \end{figure}
\end{minipage}

\end{minipage}


\medskip

Hysteresis is defined by two thresholds $\alpha<\beta$ and two values $h_1, -h_2 \in \mathbb R$ (in what follows, we are interested in the case $h_1 > 0 \geq -h_2$). Given a continuous {\em input} function $w(\tau)$, its {\em output} $\cH(w)(\tau)$ remains constant unless the input achieves the lower threshold $\alpha$ or
the upper threshold $\beta$. In the former case, the output either
{\em switches} to $h_1$ if it was equal to $-h_2$ ``just before''
or otherwise remains $h_1$. Analogously, in the latter case, the
output either {\em switches} to $-h_2$ if it was equal to $h_1$
``just before'' or otherwise remains $-h_2$. Since the function
$v(x,\tau)$ in~\eqref{eqPrototype} depends not only on $\tau$, but
also on the spatial variable $x$, one defines
$\cH(v)=\cH(v(x,\cdot))(\tau)$ ``pointwise'', i.e., for each fixed
$x$. Thus, the hysteresis operator $\cH$ becomes {\em spatially
distributed}.

Problem~\eqref{eqPrototype} is the simplest model of a
reaction-diffusion process in which a diffusive substance with
density $v(x,\tau)$ interacts in a hysteretic way with a
non-diffusive substance that affects the diffusive one via the
reaction term taking values $h_1$ or $-h_2$. The first model of such
a type was suggested by Hoppensteadt and
J{\"a}ger~\cite{HoppensteadtJaeger}. It consisted of two
reaction-diffusion equations  and one ordinary differential equation
and described the concentric rings pattern that occurs in a colony
of bacteria (Salmonella typhimurium) on a Petri plate (Fig.~\ref{figRings}).


Numerical simulations in~\cite{HoppensteadtJaeger,
HoppensteadtJaegerPoppe} yielded a pattern that was consistent with
experiments, however the rigorous mathematical description of the
model was lacking. To begin with, the well-posedness was an open
question, due to the discontinuous nature of the hysteresis
operator. First analytical results were obtained in~\cite{Alt,
VisintinSpatHyst} (see also~\cite{AikiKopfova,Kopfova, VisintinBook}
and a recent survey \cite{VisintinTen}), where existence of
solutions for multi-valued hysteresis was proved. Formal asymptotic
expansions of solutions were recently obtained in a special case
in~\cite{Ilin}. Questions about the uniqueness of solutions and
their continuous dependence on the initial data as well as a
thorough analysis of pattern formation still remained open.

In~\cite{GurTikhNonlinAnal,GurTikhSIAM13}, we formulated the
so-called {\em transversality} condition for the initial data
$\varphi(x)$ in~\eqref{eqPrototype}  that guaranteed existence,
uniqueness, and continuous dependence of solutions on initial data
for scalar equations with hysteresis. In~\cite{GurTikhMathBoh14},
this condition was generalized to systems, and
in~\cite{CurranMasterThesis} to the case $x\in\bbR^2$. For
problem~\eqref{eqPrototype}, the transversality loosely speaking
means that if $\varphi(x_0)=\alpha$ or $\varphi(x_0)=\beta$ for
some $x_0\in(-1,1)$, then $\varphi'(x_0)\ne0$. Due
to~\cite{GurTikhNonlinAnal,GurTikhSIAM13,GurTikhMathBoh14}, either the
solution  exists and is unique for all $\tau\in[0,\infty)$,
or there is $T>0$ such that the solution exists and is unique for
$\tau\in[0,T]$ and $v(x,T)$ is not transverse. The approach of~\cite{GurTikhNonlinAnal,GurTikhSIAM13,GurTikhMathBoh14} was based
on treating the problem with transverse initial data as a special  free boundary problem.
The study of regularity of the
emerging free boundary was initiated in~\cite{AU1, AU2}. For an overview on
classical free boundary problems of both elliptic and parabolic types,
we refer the reader to \cite{CaffarelliBook, SUW, UraltsevaFreeBoundaryBook} and the
references therein.

The key question which we address in this paper is how the solution
may behave after it becomes nontransverse. To answer this question,
we consider the nontransverse initial data. First, set $\beta=0$ (without loss of generality) and consider an initial function  $\varphi(x)=-cx^2+o(x^2)$ in a neighborhood $\mathcal B(0)$ of
$x=0$. By taking a smaller neighborhood if needed, we have   $\varphi(x)<0$ for $x\in \mathcal B(0)\setminus\{0\}$. We define the hysteresis at the initial moment in this
neighborhood as follows: $\cH(\varphi(0))=-h_2$ and
$\cH(\varphi(x))=h_1$ for $x\ne 0$. Now we ``regularize'' the
parabolic equation in $\mathcal B(0)$ by discretizing the spatial
variable: for any $\ep > 0$, setting
$v_n(\tau;\varepsilon):=v(\varepsilon n,\tau)$, we replace the
continuous model~\eqref{eqPrototype} in $\mathcal B(0)$ by the discrete one
\begin{equation}\label{eqDiscretePrototype}
\left\{
\begin{aligned}
& \dfrac{d v_n}{d \tau}  = \dfrac{\Delta v_n}{ \varepsilon^2}+\cH(v_n),\quad \tau>0, \ n=-N_\varepsilon,\dots,N_\varepsilon,\\
& v_n(0)=  - c(\varepsilon n)^2 + o(\varepsilon^2 n^2), \quad
n=-N_\varepsilon,\dots,N_\varepsilon,
\end{aligned}\right.
\end{equation}
 where $\Delta v_n:=v_{n+1}-2v_n+v_{n-1}$ and $N_\varepsilon\to\infty$ as
$\varepsilon\to 0$. Since we are interested in small $\varepsilon$
and in the behavior near the threshold $\beta=0$ (i.e., in a small neighborhood $\mathcal B(0)$), we consider the
next approximation by omitting $o(\varepsilon^2 n^2)$ in the initial data, replacing $N_\varepsilon$ by $\infty$, and
formally setting $\alpha:=-\infty$.
Thus,~\eqref{eqDiscretePrototype} assumes the form
\begin{equation}\label{eqDiscretePrototypeZ}
\left\{
\begin{aligned}
& \dfrac{d v_n}{d \tau}  = \dfrac{\Delta v_n}{ \varepsilon^2}+\cH(v_n),\quad \tau>0,\ n\in\bbZ,\\
& v_n(0)=  - c(\varepsilon n)^2, \quad n\in\bbZ,
\end{aligned}\right.
\end{equation}
the hysteresis operator is represented by
Fig.~\ref{figRelay}.b (see the rigorous definition in
Section~\ref{secSetting}).

A nontrivial dynamics occurs in the case  $h_1>2c>0\ge -h_2$. To
indicate the difficulty, note that, due to the initial
configuration of hysteresis, we have $\dfrac{d
v_0}{d \tau}(0;\varepsilon)=-h_2-2c<0$, but $\dfrac{d
v_n}{d \tau}(0;\varepsilon)=h_1-2c>0$ for $n\in\bbZ\setminus\{0\}$.
Thus, for small $\tau>0$, $v_0(\tau;\varepsilon)$ decreases, while
all the other {\em nodes} $v_n(\tau;\varepsilon)$,
$n\in\bbZ\setminus\{0\}$, increase. It is not clear at all, which
node achieves the threshold $\beta=0$ and switches first and hence
what a further dynamics is.

\begin{figure}[ht]
\centering
\includegraphics[width=0.90\linewidth]{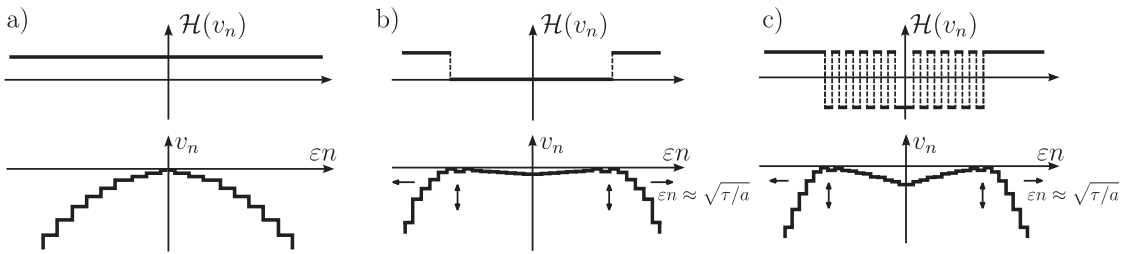}
\caption{Upper graphs represent spatial profiles of the hysteresis
$\cH(v_n)$ and lower graphs the spatial profiles of the solution
$v_n$. a) Nontransverse initial data. b) Spatial profiles at a
moment $\tau>0$ for $h_2=0$. c) Spatial profiles at a moment
$\tau>0$ for $h_2=h_1 > 0$.}
\label{figRattling}
\end{figure}

In fact, numerical analysis does not reveal any general rule that
could describe the behavior of $v_n(\tau;\varepsilon)$ for small
$\tau$. However, it reveals the formation of quite a specific
spatio-temporal pattern for large $\tau$, see
Fig.~\ref{figRattling}.


If $h_2=0$, then each node eventually achieves the threshold
$\beta=0$ and thus $\cH(v_n)$ eventually switches from $h_1$ to
$h_2=0$ for each $n\in\bbZ$. If $h_2 > 0$, then some nodes achieve the
threshold and some do not. If we denote by $N_1(j)$ and $N_2(j)$ the
number of nodes in the set $\{v_0, v_{\pm1},\dots,v_{\pm j}\}$ that
switch and do not switch, respectively, on the time interval
$[0,\infty)$, then numerics suggests that
\begin{equation}\label{eqDiscretePattern}
\lim\limits_{j\to\infty}\dfrac{N_2(j)}{N_1(j)}=\dfrac{h_2}{h_1}.
\end{equation}
Moreover, if $h_2/h_1=p_2/p_1$, where $p_1$ and $p_2$ are
co-prime integers, then, for any $j$ large enough, the set
$\{v_{j+1},\dots, v_{j+p_1+p_2}\}$ contains exactly $p_1$ nodes
that switch and $p_2$ nodes that do not switch on the time
interval $[0,\infty)$.


The next numerical observation is as follows. Let $\tau_n=\tau_n(\varepsilon)$ be the switching moment
of the node $v_n(\tau;\varepsilon)$ if this node switches  on the
time interval $[0,\infty)$ and $\tau_n:=\infty$ otherwise. Then,
for any fixed $h_2 \ge 0$, the $\tau_n$'s that are finite satisfy,
as $n\to\infty,$
\begin{equation}\label{eqTauAsymp}
\tau_n=a(\varepsilon n)^2+
\begin{cases}
\varepsilon^{2}
O(\sqrt{  n}) & \text{if } h_2=0,\\
\varepsilon^2 O(  n) & \text{if } h_2 > 0,
\end{cases}
 \end{equation}
where  $a>0$ depends on $h_1/c$ but {\em does not depend} on $h_2$
or $\varepsilon$ and $O(\cdot)$ does not depend on~$\varepsilon$.

\begin{remark}\label{remRescaling}
In Section~\ref{secStructure}, we will show that $\varepsilon$
in~\eqref{eqDiscretePrototypeZ} can be scaled out, see
scaling~\eqref{eqRescale}. In particular, all the numerical
observations concerning the dynamics of $v_n$ have been done for
$\varepsilon=1$ and then transferred to an arbitrary $\varepsilon$
according to the scaling in~\eqref{eqRescale}.
\end{remark}

Consider the function
$$
H(x,\tau;\varepsilon):=\cH(v_n(\cdot;\varepsilon))(\tau),\quad
x\in[\varepsilon n-\varepsilon/2,\varepsilon n + \varepsilon/2), \
n\in\bbZ,
$$
which is supposed to approximate the hysteresis
$\cH(v(x,\cdot))(\tau)$ in~\eqref{eqPrototype}.
Assuming the dynamics \eqref{eqDiscretePattern} and \eqref{eqTauAsymp} and taking into account
Remark \ref{remRescaling}, we see that $H(x, \tau; \ep)$
has no pointwise limit as $\varepsilon\to 0$, but converges
in a certain weak sense to the function $H(x,\tau)$ given by
$H(x,\tau)=0$ for $\tau>ax^2$ and $H(x,\tau)=h_1$ for $\tau<ax^2$.
We emphasize that $H(x,\tau)$ does not
depend on $h_2$ (because $a$ does not). On the other hand, if
$h_2 > 0$, the hysteresis operator $\cH(v(x,\cdot))(\tau)$
in~\eqref{eqPrototype} cannot take value $0$ by definition, which
clarifies the essential difficulty with the well-posedness of the
original problem~\eqref{eqPrototype} in the nontransverse case. To
overcome the non-wellposedness, one need to allow the intermediate
value $0$ for the hysteresis operator.

Such a re-definition of hysteresis is consistent with the behavior
of $v(x,\tau)$ (also observed numerically) in the following sense. For a fixed
$\varepsilon>0$, the spatial profile of
$v_n(\cdot;\varepsilon)(\tau)$ forms two humps propagating away
from the origin according to~\eqref{eqTauAsymp}. The cavity
between the humps has a bounded steepness characterized by the
relations
 \begin{equation}\label{eqGradient}
|v_{k+1}(\tau;\varepsilon)-v_{k}(\tau;\varepsilon)|\le b
\varepsilon^2,\quad |k|\le n,\ \tau\ge \tau_n,\ n=0,1,2,\dots,
\end{equation}
where $b>0$ does not depend on $k$, $n$, and $\varepsilon$. As
time goes on, the profile executes downwards and upwards motions,
always remaining beneath the threshold $\beta=0$ and  hitting this
threshold at specific nodes characterized
by~\eqref{eqDiscretePattern}. We call such a behavior of $v_n$ and
$\cH(v_n)$ {\em rattling}. Furthermore, numerics indicates that,
as $\varepsilon\to 0$, the function
$$
V(x,\tau;\varepsilon):=v_n(\tau;\varepsilon),\quad
 x\in[\varepsilon n-\varepsilon/2,\varepsilon n + \varepsilon/2),\ n\in\bbZ,
$$
approximates a smooth function $V(x,\tau)$, which satisfies
$V(x,\tau)=0$ for $\tau>ax^2$ due to~\eqref{eqTauAsymp}
and~\eqref{eqGradient}. In other words, $V(x,\tau)$ sticks to the
threshold line $\beta=0$ on the expanding interval
$x\in(-\sqrt{\tau/a},\sqrt{\tau/a})$.

We recall paper~\cite{Alt}, in which Alt proved the existence of a
function $V(x,\tau)$ that satisfies the equation
$$
V_\tau = V_{xx} + \gamma(x,\tau),
$$
where $\gamma(x,\tau) = \cH(V(x,\cdot))(\tau)$ a.e. on the set $A:=\{(x,\tau):
V(x,\tau)\ne\alpha,\beta\}$ and $\gamma(x,\tau)=0$ a.e. on the
set $B:=\{(x,\tau): V(x,\tau)=\alpha\ \text{or}\ \beta\}$ (which potentially may have a nonzero measure). Thus, our
heuristic argument provides a qualitative description of the sets
$A$ and $B$ and justifies the completion of hysteresis by the zero
value via the thermodynamical limit. To make this argument
mathematically rigorous, we should first rigorously describe the
rattling phenomenon in the discrete
system~\eqref{eqDiscretePrototypeZ}. This is the central topic of
the present paper, in which we concentrate on the case $h_2=0$ and
develop general tools for treating discrete reaction-diffusion
equations with discontinuous hysteresis. The application of these
tools to the case $h_2 > 0$ will be a subject of a forthcoming paper. We expect that these tools will be applicable whenever $h_2/h_1$ is rational.

\begin{figure}[ht]
{\ \hfill\epsfxsize120mm\epsfbox{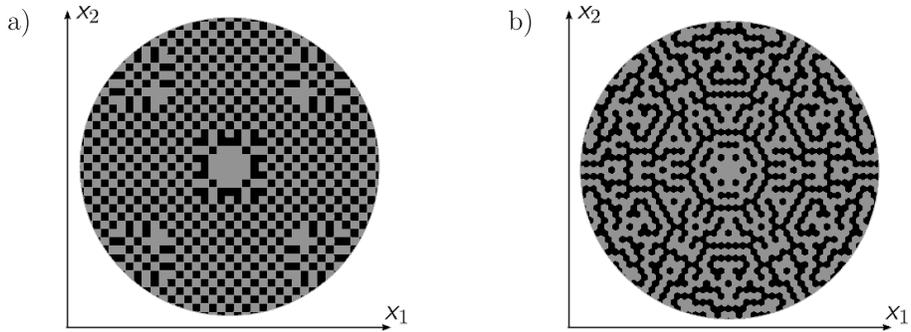}\hfill\ }
\caption{A snapshot for a time moment $\tau>0$ of a
two-dimensional spatial profile of hysteresis taking values
$h_1>4c>0$ and $h_2=h_1 > 0$. The nontrasverse initial data is
given by $\varphi(x)=-c(x_1^2+x_2^2)$. Grey (black) squares or
hexagons correspond to the nodes that have (not) switched on the
time interval $[0,\tau]$. a) Discretization on the square lattice.
b) Discretization on the triangular lattice.} \label{figMultiDim}
\end{figure}

Before we proceed with the description of our tools and of the
structure of the paper, let us make two more comments. First, the
rattling phenomenon also occurs in multidimensional domains. For
example, Fig.~\ref{figMultiDim} illustrates the switching pattern
for a two-dimensional analog of~\eqref{eqDiscretePrototypeZ},
where we have implemented spatial discretizations on the square
and triangular lattices, respectively.
Moreover,   numerical analysis of the Hoppensteadt--J{\"a}ger system
indicates that the solution remains transverse as long as the
central disc in Fig.~\ref{figRings} gets formed, but the formation
of all the rings occurs via rattling.

Second, the rattling phenomenon is not a pure consequence of a
discontinuous nature of hysteresis, but rather a consequence of
bistability in a system. In particular, it persists in bistable
slow-fast reaction-diffusion systems. The simplest example is the
system
\begin{equation}\label{eqPrototypeSlowFast}
 v_\tau  = v_{xx}+w,\qquad
 \delta w_\tau=f(v,w),
\end{equation}
where $\delta>0$ is a small parameter and the nullcline of $f(v,w)$
is $Z$- or $S$-shaped. Formally, system~\eqref{eqPrototypeSlowFast}
can be treated as another regularization of
system~\eqref{eqPrototype}. In the case where the nullcline of
$f(v,w)$ is $S$-shaped, one should replace $h_1$ and $-h_2$ in the
definition of hysteresis $\cH(v)$ by appropriate functions $H_1(v)$
and $H_2(v)$, see Fig.~\ref{figSShapedf}.


As $\delta\to 0$, the spatial profiles of $v$ and $w$
in~\eqref{eqPrototypeSlowFast} behave similarly to
$V(x,\tau;\varepsilon)$ and $H(x,\tau;\varepsilon)$, respectively,
as $\ep \to 0$, see Fig.~\ref{figSlowFastPattern}, with the
exception that the profile of $w$ remains continuous and forms steep
transition layers between mildly sloping steps of width tending to
$0$ as $\delta \to 0$.
Interestingly, the time-scale separation parameter $\delta$
in~\eqref{eqPrototypeSlowFast} yields the same effect as the
grid-size parameter $\ep$ in~\eqref{eqDiscretePrototypeZ}. As far
as we know, such a rattling phenomenon for slow-fast systems has not
been explained in the literature, either.

\noindent
\begin{minipage}{\linewidth}
      \centering
\begin{minipage}{0.49\linewidth}
    \begin{figure}[H]
    {\ \hfill\epsfxsize81mm\epsfbox{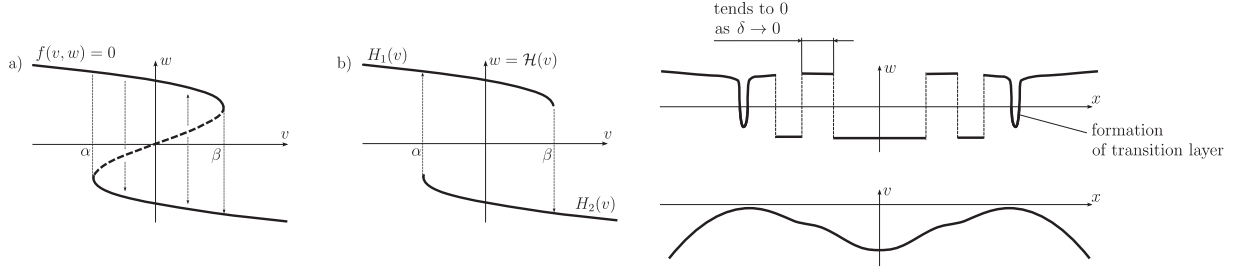}\hfill\ }  \caption{a) The nullcline of the $S$-shaped nonlinearity $f(v,w)$.
    b) Hysteresis with nonconstant branches $H_1(v)$ and $H_2(v)$.} \label{figSShapedf}
    \end{figure}
\end{minipage}
      \hspace{0.01\linewidth}
\begin{minipage}{0.48\linewidth}
    \begin{figure}[H]
    {\ \hfill\epsfxsize75mm\epsfbox{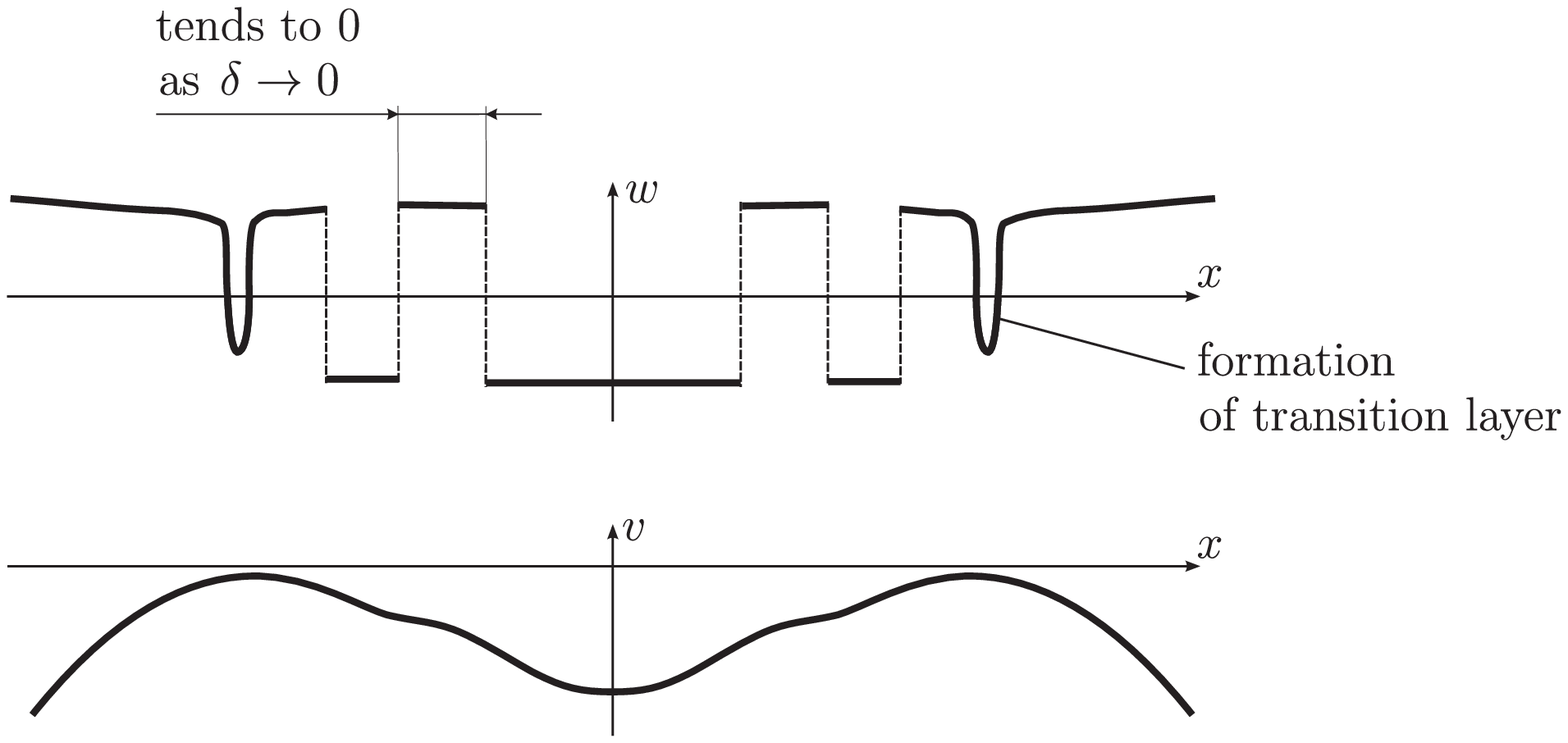}\hfill\ }
    \caption{Upper and lower graphs are spatial profiles of the solution
    $w(x,\tau)$ and $v(x,\tau)$, respectively, for
    problem~\eqref{eqPrototypeSlowFast} with initial data
    $v|_{\tau=0}=-cx^2 + o(x^2)$, $w|_{\tau=0}=h_1$.} \label{figSlowFastPattern}
    \end{figure}
\end{minipage}
\end{minipage}

\subsection{Structure of the paper}\label{secStructure}

Now we come back to the main topic of this paper, namely, discrete
system~\eqref{eqDiscretePrototypeZ}. As it was
mentioned in Remark~\ref{remRescaling}, $\varepsilon$
in~\eqref{eqDiscretePrototypeZ} can be scaled out. Indeed, setting
\begin{equation}\label{eqRescale}
 t :=
\varepsilon^{-2} \tau,\quad   u_n(t):=\varepsilon^{-2}
v_n(\tau;\varepsilon)
\end{equation}
and using the equalities (recall that $\alpha = -\infty$ and
$\beta=0$)
$$
\cH(v_n)(\tau) = \cH(\varepsilon^2
u_n(\varepsilon^{-2}\cdot))(\tau) =
\cH(u_n(\varepsilon^{-2}\cdot))(\tau) =
\cH(u_n)(\varepsilon^{-2}\tau) = \cH(u_n)(t),
$$
 we can
rewrite~\eqref{eqDiscretePrototypeZ} as follows:
\begin{equation}\label{eqDiscretePrototypeZEps1}
\left\{
\begin{aligned}
& \dot u_n  = \Delta   u_n+\cH(  u_n),\quad t>0,\ n\in\bbZ,\\
&   u_n(0)=  - c n^2, \quad n\in\bbZ,
\end{aligned}\right.
\end{equation}
where $\dot{}=d/d t$. Problem \eqref{eqDiscretePrototypeZEps1} does not involve $\ep$,
which justifies the fact that $u_n(t)$
in~\eqref{eqRescale} does not depend on $\varepsilon$. Note that $c$
in~\eqref{eqDiscretePrototypeZEps1} could be also scaled out replacing $u_n(t)$, $h_1$ and $-h_2$ by $c\tilde u_n(t)$, $c\tilde h_1$ and $-c\tilde h_2$, respectively.
We prefer not to do this, in order to keep track of what
exactly is influenced in our intermediate calculations by the
``tangency'' constant~$c$.

From now on, we concentrate on the case $h_2=0$. Due
to~\eqref{eqTauAsymp} and~\eqref{eqRescale}, the asymptotics for
the switching moment $t_n$ of the node $u_n(t)$ is expected to be
\begin{equation}\label{eqtnExpected}
t_n=a n^2 + q_n,\quad |q_n|\le E\sqrt{n},
\end{equation}
where $E>0$ does not depend on $n\in\bbZ$.


Our {\em main result} (Theorem~\ref{thMainResult}) is as follows. Let $h_1>2c>0$ and
$h_2=0$. Assume that
\begin{equation}\label{eqFinitelyManySwitchings}
\begin{aligned}
 &\text{finitely many nodes}\ u_n(t)\\
 & \text{switch at  moments}\ t_{n},\ n=0,1,\dots,n_0,\
\text{satisfying~\eqref{eqtnExpected}},
\end{aligned}
\end{equation}
 where the constants
$a=a(h_1/c)>0$ and $n_0=n_0(E)=n_0(E,h_1,c)$ will be
explicitly specified in the main text. Then each node $u_n(t)$, $n\in\bbZ$,
switches; moreover, the switching occurs at a time moment $t_n$
satisfying~\eqref{eqtnExpected}.

\begin{figure}[ht]
\centering
\includegraphics[width=0.70\linewidth]{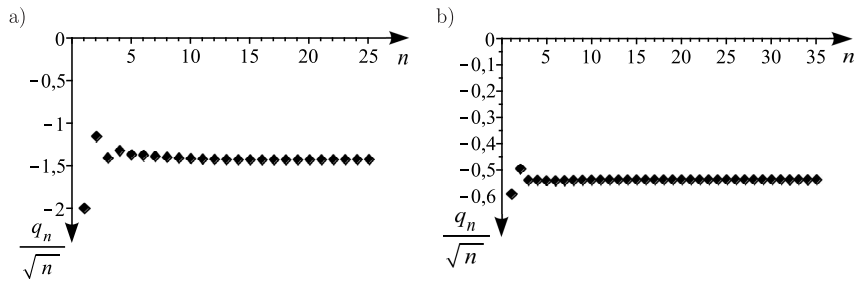}
\caption{Values of $q_n/\sqrt{n}$ for $c = 1/2$ and a) $h_1 = 1.5$, b) $h_1 = 2.0$.}
\label{figqn}
\end{figure}

Since we will provide an explicit formula for the solution $u_n(t)$,
the fulfillment of finitely many
assumptions~\eqref{eqFinitelyManySwitchings} can be verified
numerically with an arbitrary accuracy for any given values of $h_1$
and $c$ (see Fig.~\ref{figqn}).

The paper is organized as follows. In Section~\ref{secSetting}, we
give definitions for the hysteresis operator and for the solution
of problem~\eqref{eqDiscretePrototypeZEps1}. Next, we formulate
the existence and uniqueness theorem
(Theorem~\ref{thExistenceUniqueness}), which includes a
representation of the solution $u_n(t)$ via the 
\textit{discrete
Green function} $y_n(t)$. In particular,
Theorem~\ref{thExistenceUniqueness} implies that
$u_n(t)=u_{-n}(t)$, $n\in\bbZ$.

In Section~\ref{secMainResult}, we formulate our main result
(Theorem~\ref{thMainResult}).

In Section~\ref{secAux}, we formulate three main
ingredients for the proof of the main result.
\begin{enumerate}
\item The first ingredient is asymptotic formulas for the Green
function $y_n(t)$ and for its time derivatives, which were derived
in~\cite{AsympGreenFunc}.

\item The second ingredient is three equations for finding the
constant $a$ from equation~\eqref{eqtnExpected}. The equivalence of these
equations as well as the existence and uniqueness of their root
$a>0$ are proved in Appendix~\ref{secaInt}.

\item The third ingredient is the approximation of some singular
integrals by Riemann sums and corresponding error estimates, which
are proved in~\cite{ErrorEstimates2017}.
\end{enumerate}

Sections~\ref{secUntn}, \ref{secNablaUn},
and~\ref{secEstimatesUn1} are  three key steps in the proof of our
main result. The scheme of the proof is inductive. Assume we have
proved that $t_0,t_1,\dots,t_{n-1}$ satisfy~\eqref{eqtnExpected}
for some fixed $n\ge n_0+1$. We {\em fix} the hysteresis
configuration, i.e., set $H_n:=\cH(u_n)(t_{n-1})$ and consider the
solution $v_n(t)$ of the problem
$$
\left\{
\begin{aligned}
& \dot v_n=\Delta v_n + H_n,\quad t>t_{n-1},\ n\in\bbZ,\\
& v_n(t_{n-1})=u(t_{n-1}),\quad n\in\bbZ
\end{aligned}\right.
$$
(we abuse the notation by using the same letter $v$ as in
Section~\ref{subsecBackground}). Obviously, $v_n(t)=u_n(t)$ as
long as the nodes $v_n(t), v_{n+1}(t), v_{n+2}(t),\dots$ remain
below the threshold $\beta=0$.

The main theorem of Section~\ref{secUntn} (Theorem~\ref{thqn})
claims that the equation $v_n(an^2+q_n)=0$ has a root $q_n$
satisfying~\eqref{eqtnExpected}. To prove this, we use an explicit
representation of $v_n(t)$ via the convolution of $H_n$ with the
Green function $y_n(t)$ (see~\eqref{eqvInsteadOfu}). Then we use
asymptotic formulas for $y_n(t)$ (the first ingredient from
Section~\ref{secAux}) and replace the convolution by a singular
integral (the third ingredient from Section~\ref{secAux}). As
a result, we obtain a leading order term of order $n^2$, which
depends only on $a$ and $h_1/c$, and a remainder of order
$\sqrt{n}$, which also depends on $q_0,q_1,\dots,q_{n-1}$ (that are
known due to the inductive hypothesis) and on the unknown $q_n$. It
appears that the coefficient at $n^2$ vanishes due to the choice of
$a$ (the second ingredient from Section~\ref{secAux}). The
hard part is to show that the remainder vanishes for some $q_n$
satisfying~\eqref{eqtnExpected}. This is done by an application of
Brouwer's fixed-point theorem.

The time moment $t_n:=an^2+q_n$ given by Theorem \ref{thqn} is a {\em candidate} for being the
switching moment of $u_n(t)$. To show that it is the switching
moment, we have to prove that neither of the nodes $v_{n+1}(t),
v_{n+2}(t),\dots$ achieves the value $\beta=0$ on the interval $(t_{n-1},t_n]$, while $v_{n}(t)$ achieves it at the
moment $t_n$ for the first time. This is done in
Sections~\ref{secNablaUn} and~\ref{secEstimatesUn1}.

In Section~\ref{secNablaUn}, we prove that $v_{n+1}(t_n)<0$
(Theorem~\ref{thNablaun}). To do so, we estimate the gradient
$\nabla v_n(t_n):=v_{n+1}(t_n)-v_n(t_n)$ by using the representation
of $\nabla v_n(t)$ via the gradient $\nabla y_n(t)$ of the Green
function, applying asymptotic formulas for $\nabla y_n(t)$ (recall the first ingredient from Section \ref{secAux}) and again replacing the corresponding convolution by an integral (recall the third ingredient from Section \ref{secAux}). It appears that the leading order term of order $n$ vanishes  due to the second
ingredient from Section~\ref{secAux}. Thus, we calculate the
next term in the asymptotics, which turns out to be $-3h_1/4<0$.
Hence, $v_{n+1}(t_n)=\nabla v_n(t_n)\le -3h_1/8<0$.

In Section~\ref{secEstimatesUn1}, we first show that $v_n(t)$ does
not achieve the threshold $\beta=0$ for $t\in (t_{n-1},t_n)$
(Theorem~\ref{thNoSwitchBefore}). To do so, we divide the interval
$(t_{n-1},t_n)$ into two parts. We prove that the function $\dot
v_n(t)$ is so small on the first interval that it cannot overcome
the distance exceeding $-3h_1/8$ (the value coming from Theorem \ref{thNablaun} with $n+1$ replaced by $n$).
Then we prove that $\ddot
v_n(t)$ is nonnegative on the second interval. Hence, the equation
$v_n(t)=0$ has a unique root, which must be $t_n$. In particular,
$v_n(t)<0$ for $t\in (t_{n-1},t_n)$. Finally, we show that $
\nabla v_j(t)<0$ for all $t\in(t_{n-1},t_n]$ and $j\ge n$, which
implies that the nodes $v_{n+1}(t), v_{n+2}(t),\dots$ remain
negative for $t\in (t_{n-1},t_n]$
(Theorem~\ref{thNoSwitchBeforeAllun}).

In Section~\ref{secProofMainTheorem}, we combine the results from
Sections~\ref{secUntn}, \ref{secNablaUn},
and~\ref{secEstimatesUn1} and rigorously implement the inductive
scheme, which completes the proof of the main result, namely,
Theorem~\ref{thMainResult}.

The crucial role in our main result (Theorem~\ref{thMainResult})
is played by the number $n_0 = n_0(E)$, which determines the number of
switchings one has to check ``by hand'' (see \eqref{eqFinitelyManySwitchings}). The number $n_0(E)$ is
determined explicitly by 12 inequalities that must
hold for $n\ge n_0(E)$. Each inequality is referred to as a {\bf
requirement} and is introduced in the text where it is used for
the first time. These 12 requirements contain constants that are
also introduced in the text where they are used for the first
time. For reader's convenience, we have collected all those
constants in Appendices~\ref{subsecConstNoAE}--\ref{subsecConstE}
and the 12 requirements
in Appendix~\ref{refAppendixSubsectionAssumptionsN0}.

\begin{figure}[ht]
\begin{center}
      \includegraphics[width=0.95\linewidth]{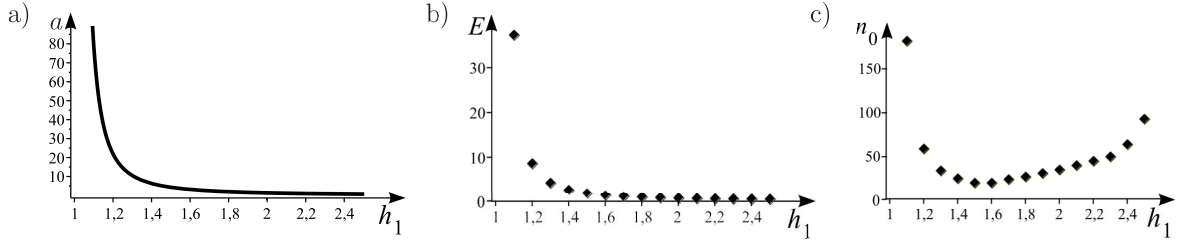}
\end{center}
\vspace{-10pt}
\caption{Dependence on $h_1$ of the values of $a$, $E$, and
$n_0(E)$ that fulfill assumptions~\eqref{eqFinitelyManySwitchings}
for $c=1/2$. a) The values of $a$ are found explicitly for all $h_1>1$ (see Section \ref{sectaInt} below). b),~c)~The values of $E$ and
$n_0(E)$ are calculated numerically for $h_1=1.1, 1.2, 1.3,  \dots, 2.5$.}
\label{figaEn0}
\end{figure}

The graphs in Fig.~\ref{figaEn0} represent the values of $a$, $E$, and
$n_0(E)$ that fulfill assumptions~\eqref{eqFinitelyManySwitchings}
for $c=1/2$ and $h_1=1.1, 1.2, 1.3,  \dots, 2.5$.


\section{Setting of the problem and a proof of its well-posedness}\label{secSetting}

For a sequence $\{v_n\}_{n\in\bbZ}$ of real numbers, we use the
notation
$$
\nabla v_n:=v_{n+1}-v_n,\qquad \Delta v_n:=\nabla v_{n}-\nabla
v_{n-1}=v_{n-1}-2v_n+v_{n+1}.
$$
Let $\{u_n(t)\}_{n\in\bbZ}$ be real-valued functions
defined for $t\ge 0$. We study the problem
\begin{align}
& \dot u_n=\Delta u_n+\cH(u_n), \quad t>0,\ n\in\bbZ, \label{equ_nEquation1}\\
& u_n(0)=-c n^2,\quad n\in\mathbb Z, \label{equ_nEquation2}
\end{align}
where $c>0$ and $\cH(w)(t)$, $t \geq 0$, is the  hysteresis operator defined for
functions $w \in C[0,\infty)$ such that $w(0)\le 0$ by
\begin{equation}\label{nonideal}
 \cH(w)(t) :=
 \begin{cases} h_1 & \text{if } w(s)< 0 \ \text{for all } s \in [0,t], \\
              0 & \text{if }  w(s)=0 \ \text{for some}\  s \in  [0,t],
 \end{cases}
\end{equation}
where $h_1>0$ is fixed. In other words, the  output of hysteresis is $h_1$
unless the input achieves the zero threshold; at this moment, the
hysteresis {\it switches} and since then the  output of hysteresis
remains $0$. In the context of
problem~\eqref{equ_nEquation1}--\eqref{nonideal}, we will
say ``a node $u_n(t)$ switches'' or ``a node $n$
switches'' whenever $u_n(t)$ achieves the value zero for the first
time.

\begin{remark}
In the terminology of, e.g.,~\cite{KrasnBook,VisintinBook}, the
hysteresis operator~\eqref{nonideal} is a non-ideal relay with the
thresholds $-\infty$ and $0$; see Fig.~\ref{figRelay}.b.
\end{remark}

From now on, we assume throughout that the following condition holds.

\begin{condition}\label{condHc}
$h_1>2c>0$.
\end{condition}

We note that the function $\cH(v)(t)$ may have discontinuity
(actually, at most one) even if $v\in C^\infty[0,\infty)$.
Therefore, one cannot expect that a solution of
problem~\eqref{equ_nEquation1}--\eqref{nonideal} is
continuously differentiable on $[0,\infty)$. Thus, we define a
solution as follows.

\begin{definition}\label{defSolution}
  We say that a sequence
$\{u_n(t)\}_{n\in\bbZ}$ is a {\em solution of
problem~\eqref{equ_nEquation1}--\eqref{nonideal} on the time
interval $(0,T)$}, $T>0$,  if
\begin{enumerate}
\item $u_n\in C[0,T]$ for all $n\in\bbZ$,

\item\label{defSolution2} for each $t\in[0,T]$, there exists $A,\alpha\ge 0$ such that
$\sup\limits_{s\in[0,t]}|u_n(s)|\le A e^{\alpha |n|}$ for all
$n\in\bbZ$,

\item there is a finite  sequence $
0=\tau_0<\tau_1<\dots<\tau_J=T$, $J\ge1$,  such that $u_n\in
C^1(\tau_j,\tau_{j+1})$ for all $n\in\bbZ$ and $j=0,\dots,J-1$,

\item the equations in~\eqref{equ_nEquation1} hold in
$(\tau_j,\tau_{j+1})$ for all $n\in\bbZ$ and  $j=0,\dots,J-1$,

\item $u_n(0)=-cn^2$ for all $n\in\bbZ$.
\end{enumerate}
  We say that a sequence
$\{u_n(t)\}_{n\in\bbZ}$ is a {\em solution of
problem~\eqref{equ_nEquation1}--\eqref{nonideal} on the time
interval $(0,\infty)$} if it is a solution on $(0,T)$ for all
$T>0$.
\end{definition}

\begin{remark}
If $\{u_n(t)\}_{n\in\bbZ}$ is a solution, then, as we have
mentioned above, the function $\cH(u_n)(t)$ has at most one
discontinuity point for each fixed $n\in\bbZ$. Hence, the
equations in~\eqref{equ_nEquation1} imply that each function $\dot
u_n(t)$ has at most one discontinuity point on $[0,\infty)$.
\end{remark}

Before we treat existence and uniqueness of a solution, let us
introduce one of our main tools, namely, the so-called {\it
discrete Green function}
\begin{equation}\label{eqynIntegralFormula}
y_n(t):=\dfrac{1}{2\pi}\int_{-\pi}^\pi\dfrac{1-e^{-2t(1-\cos\theta)}}{2(1-\cos\theta)}\,e^{i
n\theta}\,d\theta, \quad t\ge 0.
\end{equation}
One can directly check that $y_n\in C^\infty[0,\infty)$ and
$y_n(t)$ solves the problem
\begin{equation}\label{eqy_nGreen}
\left\{
\begin{aligned}
& \dot y_0=\Delta y_0+1, & & t>0,\\
& \dot y_n=\Delta y_n, & & t>0,\ n\ne0,\\
& y_n(0)=0, & & n\in\mathbb Z.
\end{aligned}
\right.
\end{equation}

Below, we will use the fact that
\begin{equation}\label{eqDotynMonotoneInN}
\dot y_{n+1}(t)<\dot y_n(t),\quad t>0,\ n=0,1,2,\dots,
\end{equation}
which follows from the formula $\dot y_n(t)=e^{-2t}I_n(2t)$, where $I_n(s)$ is the modified Bessel function of the first kind (see~\cite[Sec.~9.6.19]{Abramowitz}), and from, e.g.,~\cite{Amos}. We will also use the estimate, which follows from the series
representation of the modified Bessel
function~\cite[Sec.~9.6.10]{Abramowitz}:
\begin{equation}\label{eqEstimateynSeries}
0\le \dot y_n(t)=e^{-2t}I_n(2t)=e^{-2t}
t^{|n|}\sum\limits_{m=0}^\infty \dfrac{t^{2m}}{m!(m+|n|)!}\le
\dfrac{e^{t^2-2t} t^{|n|}}{|n|!},\quad n\in\bbZ.
\end{equation}

Below we prove the following existence and uniqueness result.

\begin{theorem}\label{thExistenceUniqueness}
\begin{enumerate}
 \item Problem~\eqref{equ_nEquation1}--\eqref{nonideal}
has a unique solution $\{u_n(t)\}_{n\in\bbZ}$ on the time interval
$(0,\infty)$.

\item Let $t_n$ be the switching moment of the node $u_n(t)$ if
this node switches  on the time interval $[0,\infty)$ and
$t_n:=\infty$ otherwise. Then
\begin{equation}\label{eqtntoinfty}
t_n\ge \dfrac{c n^2}{h_1-2c},\quad n\in\bbZ.
\end{equation}

\item Let $S(t)$ be the set of nodes that   switch on the time
interval $[0,t]$, i.e.,
\begin{equation}\label{eqSt}
S(t):=\{k\in\bbZ: \cH(u_k)(t)=0\},
\end{equation}
and let $|S(t)|$ be the number of elements in $S(t)$. Then $S(t)$
is finite for each $t>0$, symmetric with respect to the origin,
$|S(t)|\to\infty$ as $t\to\infty$, and
\begin{equation}\label{eqSolu}
u_n(t)=-cn^2+(h_1-2c)t-h_1\sum\limits_{k\in
S(t)}y_{n-k}(t-t_k),\quad t\in[0,\infty),
\end{equation}
where we put $y_{n-k}(t-t_k)=0$ for $t<t_k$,

\item for each $n\in\bbN$, we have $t_{-n}=t_n$  and
$u_{-n}(t)\equiv u_n(t)$.
\end{enumerate}
\end{theorem}
\proof

{\bf Step 1.} Using~\eqref{eqy_nGreen}, we see that the
functions
\begin{equation}\label{eqSolTau1}
z_n^{(1)}(t):=-c n^2+ (h_1-2c)t-h_1 y_n(t),\quad n\in\bbZ,
\end{equation}
satisfy the initial condition~\eqref{equ_nEquation2} and the equation in~\eqref{equ_nEquation1} as long as
$z_n^{(1)}(t)<0$ for all $n\in\bbZ\setminus\{0\}$.
By comparing $z_n^{(1)}(t)$ with the solution
\begin{equation}\label{eqznplus}
z_n^{+}(t)=-cn^2+(h_1-2c)t
\end{equation}
of problem~\eqref{equ_nEquation1}--\eqref{nonideal} with
$\cH(u_n)$ replaced by $h_1$ for all $n\in\bbZ$, it is not
difficult to see that
\begin{equation}\label{eqzn1lessznplus}
z_n^{(1)}(t)\le z_n^{+}(t),\quad t\ge 0,\ n\in\bbZ.
\end{equation}
 Therefore, the
   time moment  $t_{n}$ at which $z_{n}(t)$ vanishes for the first time is not less than the moment $t_n^+=c n^2/(h_1-2c)$ at which
$z_n^{+}(t)$ vanishes.

In particular,  $z_n^{(1)}(t)<0$ for all $n\in\bbZ\setminus\{0\}$
and $t\in [0,c/(h_1-2c))$. Let $\tau_1:=\sup\{t>0:z_n^{(1)}(t)<0\
\forall n\in \bbZ\setminus\{0\}\}$. It follows
from~\eqref{eqSolTau1}, \eqref{eqr01t}, and~\eqref{eqA} that
$\tau_1$ is finite. As we have seen, $\tau_1\ge c/(h_1-2c)>0$.
Furthermore,  $z_n^{(1)}(t)$ satisfy the growth condition from
item~\ref{defSolution2} of Definition~\ref{defSolution} for
$t\in[0,\tau_1]$. This follows from~\eqref{eqSolTau1} and the fact
that $|y_n(t)|$ given by~\eqref{eqynIntegralFormula} are bounded
on any finite time interval, uniformly with respect to $n\in\bbZ$.
Thus, $u_n(t):=z_n^{(1)}(t)$ is a solution of
problem~\eqref{equ_nEquation1}--\eqref{nonideal} on the time
interval $(0,\tau_1)$.

Let us  prove that  the solution $u_n(t)$ is unique on
$(0,\tau_1)$. Assume we have another solution $\tilde u_n(t)$ on a
time interval $(0,\tilde\tau_1)$, where $\tilde\tau_1\le\tau_1$ is
such that $\tilde u_n(t)<0$ for all $t\in(0,\tau_1)$ and
$n\in\bbZ\setminus\{0\}$. Then the difference
$w_n(t):=u_n(t)-\tilde u_n(t)$ must satisfy the homogeneous
diffusion equation on the time interval $(0,\tilde\tau_1)$ with
the zero initial data
$$
\begin{aligned}
\dot w_n(t)&=\Delta w_n(t),\quad
t\in(0,\tilde\tau_1),\ n\in\bbZ,\\
w_n(0)&=0,\quad n\in\bbZ.
\end{aligned}
$$
If we looked for solutions that are square summable with respect to $n\in\bbZ$, then the application of the discrete  Fourier transform would immediately imply that all $w_n(t)\equiv 0$. However, we are interested in solutions that may have exponential growth with respect to $n\in\bbN$ (see item~\ref{defSolution2} in Definition~\ref{defSolution}). We will argue as follows. For each  $N\in\bbN$, we consider the functions
\begin{equation}\label{eqvwN}
\zeta_n(t)=\zeta_n^N(t):=w_n(t)\ \text{for } |n|\le N,\qquad \zeta_n(t)=\zeta_n^N(t):=0\
\text{for } |n|\ge N+1.
\end{equation}
They
satisfy the relations
\begin{equation}\label{eqvParabolic0}
\begin{aligned}
\dot \zeta_n(t)&=\Delta \zeta_n(t)+G_n^N(t),\quad
t\in(0,\tilde\tau_1),\ n\in\bbZ,\\
\zeta_n(0)&=0,\quad n\in\bbZ,
\end{aligned}
\end{equation}
  where $G_n^N(t)=0$ for $|n|\le N-1$ and $|n|\ge
N+2$, $G_{\pm N}^N(t)=w_{\pm(N+1)}(t)$, and $G_{\pm
(N+1)}^N(t)=-w_{\pm N}(t)$.
Since no more than finitely many elements in the sequences
$\{\zeta_n(t)\}_{n\in\bbZ}$ and $\{G_n^N(t)\}_{n\in\bbZ}$ are nonzero,
we can apply the discrete Fourier transform
to~\eqref{eqvParabolic0} and obtain
\begin{equation}\label{eqSolutionv}
\zeta_n^N(t)=\sum\limits_{|k|=N}^{N+1} \int_{0}^t \dot
y_{n-k}(t-s)\,G_k^N(s)\,ds,\quad t\in[0,\tilde\tau_1],\ n\in\bbZ.
\end{equation}

Now let us fix $n\ge 0$ and $t\in[0,\tilde\tau_1]$. By assumption,
there exist $A,\alpha\ge 0$ such that
\begin{equation}\label{eqwExpBounded}
\sup\limits_{s\in[0,t]}|w_k(s)|\le  A e^{\alpha |k|},\quad
k\in\bbZ.
\end{equation}

Combining~\eqref{eqvwN}, \eqref{eqSolutionv},
\eqref{eqEstimateynSeries}, \eqref{eqwExpBounded} and choosing
$N\ge n$, we have
$$
|w_n(t)|=|\zeta_n^N(t)|\le t
\sum\limits_{|k|=N}^{N+1}\sup\limits_{s\in[0,t]}\left(|\dot
y_{n-k}(s)|\cdot |G_k^N(s)|\right)\le c \dfrac{t^{N} e^{\alpha
N}}{(N-n)!}\to 0\quad \text{as } N\to\infty,
$$
where $c=c(n,t)\ge 0$ does not depend on $N$. Therefore,
$w_{-n}(t)\equiv w_n(t)\equiv 0$. This proves that $u_n(t)$ is a
unique solution of problem~\eqref{equ_nEquation1}--\eqref{nonideal} on the time interval $(0,\tau_1)$.

{\bf Step 2.} Set $S(\tau_1):=\{0\}\cup\{n\in\bbZ :
z_n^{(1)}(\tau_1)=0$\}, cf.~\eqref{eqSt}. Due to~\eqref{eqznplus}
and~\eqref{eqzn1lessznplus}, the set $S(\tau_1)$ is finite. Due
to~\eqref{eqSolTau1} and the symmetry $y_n(t)\equiv y_{-n}(t)$,
the set $S(\tau_1)$ is symmetric with respect to the origin. Note
that $t_0=0$ is the switching moment of the node $0$, while
$t_n=\tau_1$  are the switching moments of the nodes $n\in
S(\tau_1)\setminus\{0\}$.

Using~\eqref{eqy_nGreen} and assuming $y_n(t):=0$ for $t<0$, we
see that the functions
\begin{equation}\label{eqSolTau2}
\begin{aligned}
z_n^{(2)}(t)&=-cn^2+(h_1-2c)t-h_1 \left(y_n(t)+\sum\limits_{k\in
S(\tau_1)\setminus\{0\}}y_{n-k}(t-\tau_1)\right)\\
&=-cn^2+(h_1-2c)t-h_1\sum\limits_{k\in
S(\tau_1)}y_{n-k}(t-t_k),\quad t\in[0,\infty),
\end{aligned}
\end{equation}
satisfy the equations in~\eqref{equ_nEquation1} as long as
$z_n^{(2)}(t)<0$ for all $n\in \bbZ\setminus S(\tau_1)$, i.e., as
long as $S(t)=S(\tau_1)$. Obviously, $z_n^{(2)}(t)$ also satisfy
the initial condition~\eqref{equ_nEquation2}.

As in Step 1, we see that the time moment  $t_{n}$ at which
$z_{n}(t)$, $n\in \bbZ\setminus S(\tau_1)$, vanishes for the first
time is not less than $c n^2/(h_1-2c)$. Hence, there is a positive
time interval (of length bigger than $\tau_1$) on which
$z_n^{(2)}(t)<0$ for all $n\in \bbZ\setminus S(\tau_1)$.

Let $\tau_2:=\sup\{t>0:z_n^{(2)}(t)<0\ \forall n\in \bbZ\setminus
S(\tau_1)\}$. It follows from~\eqref{eqSolTau1}, \eqref{eqr01t},
and~\eqref{eqA} that $\tau_2$ is finite. We have proved that
$\tau_2>\tau_1$. Furthermore, $z_n^{(2)}(t)$ satisfy the growth
condition from item~\ref{defSolution2} of
Definition~\ref{defSolution} for $t\in[0,\tau_2]$. Thus,
$u_n(t):=z_n^{(2)}(t)$ is a solution of
problem~\eqref{equ_nEquation1}--\eqref{nonideal} on the time
interval $(0,\tau_2)$. Note that $u_n(t)=z_n^{(1)}(t)$ for
$t\in[0,\tau_1]$. The uniqueness of $u_n(t)$ on the interval
$(\tau_1,\tau_2)$ can be proved similarly to Step~1.

Continuing these steps, we   obtain the desired infinite sequence
$\{\tau_j\}_{j\ge 0}$ from
Definition~\ref{defSolution}. On each step, we compare $u_n(t)$
with $z_n^+(t)$ given by~\eqref{eqznplus} and conclude that the
switching moments satisfy $t_n\ge c n^2/(h_1-2c)$. Hence,
$\tau_j\to\infty$ as $j\to\infty$.
\endproof

\section{Main result}\label{secMainResult}

We recall that Condition~\ref{condHc} is assumed to hold throughout.
Below in the text we define $a > 0$ (see Lemma \ref{lema}), $E_0 > 0$ (see \sref{eqE0}) and an increasing function $n_0 : (E_0, \infty) \to \bbN$ (see Requirements \ref{req1}--\ref{ass8} in Section~\ref{refAppendixSubsectionAssumptionsN0}).


\begin{definition}\label{condtkForuntn}
We say that a number $E \geq E_0$ is \textit{admissible} if the following holds:
\begin{enumerate}
\item each node $u_k$, $k=0,\pm1\,\dots,\pm n_0 = \pm n_0(E)$ switches at a
moment $t_k$ satisfying
\begin{equation}\label{eqEn+ScondtkForuntn}
t_k=ak^2+q_k, \quad |q_k| \leq E\sqrt{n_0},
\end{equation}
while neither of the nodes $u_{\pm (n_0+1)}, u_{\pm
(n_0+2)},\dots$  switches on the time interval $[0,t_{n_0}]$;

\item at the switching moment $t_{n_0}$, we have
\begin{equation}\label{p14s}
u_{n_0+1}(t_{n_0}) = \nabla u_{n_0}(t_{n_0})\le -\dfrac{3h_1}{8}.
\end{equation}

\end{enumerate}
\end{definition}



The main result of this paper is as follows.
If finitely many nodes $k = 0, \dots, n_0(E)$ switch at time moments $t_k$ satisfying \eqref{eqEn+ScondtkForuntn}, then \textit{all} the nodes $n \in \mathbb Z$ will switch and their switching moments will be of order $an^2$. On Fig.~\ref{figaEn0}.b, one can see the values of admissible $E$, which we found numerically for  $c=1/2$ and $h_1=1.1, 1.2, 1.3,  \dots, 2.5$. Figures~\ref{figaEn0}.a and~\ref{figaEn0}.c depict corresponding values of $a$ and $n_0(E)$, respectively.

The rigourous formulation of our main result is as follows.

\begin{theorem}\label{thMainResult}
Assume that $E \geq E_0$ is an admissible number, and let $n_0 =
n_0(E)$. Then for all $n\ge n_0+1$$:$
\begin{enumerate}
\item Each of the nodes $u_n$  switches at a moment $t_n$
satisfying
\begin{equation}\label{p4s}
t_n=an^2+q_n, \quad |q_n| \leq E\sqrt{n},
\end{equation}
$$
t_{k}<t_{n_0}<t_{n_0+1}<\dots, \quad k = 0, 1, \dots, n_0 - 1,
$$

\item  There exists $A_{\nabla}>0$ depending on $h_1,c,E$, but not on $n$,  such that
\begin{equation*}
  \left|\nabla u_{n}(t_{n})+\dfrac{3h_1}{4}\right| \le A_{\nabla} n^{-1/2}, \qquad \nabla u_{n}(t_{n}) \le -\dfrac{3h_1}{8}.
\end{equation*}
\end{enumerate}
\end{theorem}


\section{Auxiliary Statements}\label{secAux}

In this section, we formulate several auxiliary statements. Each of them is a key ingredient in the proof of our main result, i.e., Theorem \ref{thMainResult}.

In Section \ref{sect5yn} (Proposition \ref{theoremyn}), we establish asymptotic formulas for the discrete Green function $y_n(t)$ given by \eqref{eqynIntegralFormula}. It is essential that the leading order terms in the asymptotics depend only on $n/\sqrt{t}$, while the remainders are estimated uniformly with respect to $n$.

In Section \ref{sectaInt}, we consider three expressions containing integrals \eqref{eqIntfgh} of leading order terms in the asymptotics of $y_n(t)$, $\nabla y_n(t)$, and $\dot{y}_n(t)$, respectively. These three expressions will enter the leading order terms in asymptotic formulas for $u_n(t_n)$, $\nabla u_n(t_n)$, and $\dot{u}_n(t_n)$. In Proposition \ref{lema}, we show that these terms vanish for the same value of $a$, thus determining the ``propagation rate'' $an^2$ in the switching moment asymptotics for $t_n$ in \eqref{eqEn+ScondtkForuntn} and \eqref{p4s}.

In Section \ref{sectIntApprox}, we elaborate on properties of integrals \eqref{eqIntfgh} from Section \ref{sectaInt}. In the proof of our main result,  these integrals will play the role of  approximation of some Riemann sums. Note that the corresponding integrands are not smooth functions, but have singularities of order $(1-x)^{1/2}$ or $(1-x)^{-1/2}$ at $x = 1$. In Propositions \ref{lemIntSqrt1} and \ref{lemIntSqrt2}, we provide error estimates for approximation of such integrals by their Riemann sums.

\subsection{Properties of the discrete Green function $y_n(t)$}\label{sect5yn}\label{AsympGreenFunc}

Consider the functions $h, f, g, \tilde{f}: \bbR^+ \to \bbR$ given by
\begin{equation}
\label{eqfghft}
h(x)  := \dfrac{1}{2\sqrt{\pi}}\,e^{-\frac{x^2}{4}},  \quad
f(x)  := 2x\int\limits_x^\infty y^{-2}h(y)\,dy, \quad
g(x)  := f'(x),  \quad
\tilde{f}(x)  := -\frac{h'''(x)}{6x}.
\end{equation}
Note that these functions belong to $C^{\infty}[0, \infty)$ and decay to zero as $x \to \infty$, together with all their derivatives, faster than any exponential. Moreover,
\begin{equation}\label{eqfgh}
h(x) = g'(x) = f''(x), \quad 2h(x) + xg(x) - f(x) = 0, \quad g(0) = -\frac{1}{2}.
\end{equation}
Consider the functions $r_0, \tilde{r}_1, r_1, r_2, w_0, w_1: (\bbN \cup \{0\}) \times \bbR_+ \to \bbR$ satisfying the following relations for $n = 0, 1, 2, \dots$ and $t > 0$:
\begingroup
\allowdisplaybreaks
\begin{align}
\label{eqr01t} y_n(t) & = \sqrt{t}f\lr{\dfrac{n}{\sqrt{t}}} + r_0(n, t), &
 y_n(t) & = \sqrt{t}f\lr{\dfrac{n}{\sqrt{t}}} + \frac{1}{\sqrt{t}}\tilde{f}\lr{\dfrac{n}{\sqrt{t}}} + \tilde{r}_1(n, t),\\
\label{eqr12} \dot{y}_n(t) & = \frac{1}{\sqrt{t}}h\lr{\dfrac{n}{\sqrt{t}}} + r_1(n, t), &
\ddot{y}_n(t) & = \frac{1}{t\sqrt{t}}h''\lr{\dfrac{n}{\sqrt{t}}} + r_2(n, t),\\
\nabla \dot{y}_n(t) & = \frac{1}{t} h'\lr{\dfrac{n}{\sqrt{t}}}+ w_1(n, t), &
\label{eqw01} \nabla y_n(t) & = g\lr{\dfrac{n}{\sqrt{t}}} + \frac{1}{2\sqrt{t}}h\lr{\dfrac{n}{\sqrt{t}}} + w_0(n, t).
\end{align}
\endgroup

We fix throughout the paper
\begin{equation}\label{eqtau0}
\tau_0 > 0.
\end{equation}

The following estimates are proved in \cite{AsympGreenFunc}.
\begin{prp}\label{theoremyn}
There exist constants $A_0, A_1, A_2, \tilde{A}_1, B_0, B_1, A_2^*,
B_2^* > 0$ $($depending on $\tau_0$$)$ such that, for all $t \geq
\tau_0$, $n = 0, 1, 2, \dots, $ and $i = 0, 1, 2$, the following
inequalities hold\str{:}
\begingroup
\allowdisplaybreaks
\begin{align}
\label{eqA} |r_i(n, t)| & \leq A_i \frac{1}{t^i\sqrt{t}},  & |\tilde{r}_1(n, t)| & \leq \tilde{A}_1 \frac{1}{t\sqrt{t}},  \\
\label{eqB01} |w_0(n, t)| & \leq B_0 \frac{1}{t},  &|w_1(n, t)|  & \leq B_1 \frac{1}{t\sqrt{t}}, \\
\label{eqAB2*} |\ddot{y}_n(t)| & \leq A_2^* \frac{1}{t\sqrt{t}},  &|\nabla \ddot{y}_n(t)| & \leq B_2^* \frac{1}{t^2}.
\end{align}
\endgroup
\end{prp}

\subsection{Equivalence of some equations}\label{sectaInt}
Consider the functions
\begingroup
\allowdisplaybreaks
\begin{align}
\label{eqFG} F(a, x) & := \sqrt{a(1-x^2)}\,f\left(\dfrac{1}{\sqrt{a}}\sqrt{\dfrac{1-x}{1+x}}\right), & 
 G(a, x) & := g\left(\dfrac{1}{\sqrt{a}}\sqrt{\dfrac{1-x}{1+x}}\right),\\
\label{eqHH1} H(a, x) & := \dfrac{1}{\sqrt{a(1-x^2)}}\,h\left(\dfrac{1}{\sqrt{a}}\sqrt{\dfrac{1-x}{1+x}}\right), & 
 H_1(a, x) & := \dfrac{1}{a(1-x^2)}\!h'\!\left(\dfrac{1}{\sqrt{a}}\sqrt{\dfrac{1-x}{1+x}}\right),
\end{align}
\endgroup
where $a > 0$, $x \in (-1, 1)$ and $f,g,h$ are given by~\eqref{eqfghft}.
Set
\begin{equation}
\label{eqIntfgh}
I_F(a) := \int_{-1}^1 F(a, x) dx, \quad
I_G(a) := \int_{-1}^1 G(a, x) dx, \quad
I_H(a) := \int_{-1}^1 H(a, x) dx.
\end{equation}

The following proposition is proved in Appendix \ref{secaInt}.
\begin{prp}\label{lema}
Each of the three equations
\begingroup
\allowdisplaybreaks
\begin{align}
-c+(h_1-2c)a-h_1I_F(a)&=0,\label{eqah1f}\\
-2c-h_1I_G(a)&=0,\label{eqah1g}\\
(h_1-2c)-h_1I_{H}(a)&=0\label{eqah1h}
\end{align}
\endgroup
has a unique root on the interval $(0,\infty)$. Moreover, all
these equations have the same root.
\end{prp}

In what follows, we fix $a$ given by Proposition \ref{lema} and write $F(x)$, $G(x)$, $H(x)$, $H_1(x)$, $I_F$, $I_G$, $I_H$, omitting the dependence  on $a$.

\subsection{Error estimates for Riemann sums}\label{sectIntApprox}

Let $N \in \bbN$, and let $z_0 =-1$ or $z_0 = 0$. The following propositions are proved in \cite{ErrorEstimates2017}.
\begin{prp}\label{lemIntSqrt1}
Assume that a function $F_1(x)$ can be represented as
$$
F_1(x) = c_1 (1-x)^{1/2}+ c_2 (1-x)^{3/2} + \tilde{F}_1(x),
$$
where $c_1, c_2 \in \bbR$ and $\tilde{F}_1 \in C^2[z_0, 1]$.
Denote the error estimate of the Riemann sum of the integral $\int_{z_0}^1 F_1(x) dx$ by
$$
R_n := \int_{z_0}^1 F_1(x) dx - \left( \frac{1}{2n}F_1(z_0) + \frac{1}{2n}F_1(1) + \sum_{ k = z_0 n + 1}^{n-1} \frac{1}{n}F_1\lr{\dfrac{k}{n}}\right).
$$
\begin{enumerate}
\item
There exists $L_1 = L_1(F_1, N) > 0$ such that
$$
\left|R_n\right| \leq L_1\frac{1}{n^{3/2}}, \quad n \geq N.
$$
\item If, additionally, $\tilde{F}_1(x) = c_3 (1-x)^{5/2} + c_4 (1-x)^{7/2} + \bar{F}_1(x)$, where $c_3, c_4 \in \bbR$ and $\bar{F} \in C^4[z_0, 1]$, then there exists $\bar{L}_1 = \bar{L}_1(F_1, N) > 0$ such that
$$
\left|(n+1)^2 R_{n+1} - n^2 R_n\right| \leq \bar{L}_1\frac{1}{n^{1/2}} \quad n \geq N.
$$
\end{enumerate}
\end{prp}

\begin{prp}\label{lemIntSqrt2}
Assume that a function $F_2(x)$ can be represented as
\begin{equation}\label{P7s}
F_2(x) = c_1 (1-x)^{-1/2}+ c_2 (1-x)^{1/2} + \tilde{F}_2(x),
\end{equation}
where $c_1 > 0$, $c_2 \in \bbR$, and $\tilde{F}_2 \in C^1[z_0, 1]$.
Then there exists $L_2 = L_2(F_2, N) > 0$, $L_2^* = L_2^*(F_2, N) > c_1$ and $l_2 = l_2(F_2, N)$ such that
\begin{equation}\label{eqL2}
L_2^*\frac{1}{n^{1/2}} - l_2 \frac{1}{n} \leq \int_{z_0}^1 F_2(x) dx - \sum_{k = z_0n}^{n-1} \frac{1}{n}F_2\lr{\dfrac{k}{n}} \leq
L_2\frac{1}{n^{1/2}}, \quad n \geq N.
\end{equation}
In particular, $\left|\sum \limits_{k = z_0 n}^{n-1} \frac{1}{n}F_2\lr{\frac{k}{n}}\right|$ is bounded.
\end{prp}

\begin{prp}\label{propL3}
Assume that $F_3 \in C^0[-1, 1)$ and $|(1-x)^{3/2}F_3(x)|$ is bounded on $[-1, 1)$. Then there exists $L_3 = L_3(F_3, N) > 0$ such that
$$
\left|\sum_{|k|\leq n-1} \frac{1}{n}F_3\lr{\frac{k}{n}}\right| \leq L_3{n^{1/2}}, \quad n \geq N.
$$
\end{prp}

\section{Asymptotics for $u_n(t_n)$}\label{secUntn}

\subsection{Preliminaries}

Set ($h$ is defined in \eqref{eqfghft})
\begin{equation}\label{eqfDa}
h_a(x) := h\lr{\frac{x}{\sqrt{a}}} +
h\lr{\frac{1}{x\sqrt{a}}}, \quad
D_a := \inf_{x \in (0, 1]}h_a(x)> h(a^{-1/2}),
\end{equation}
\begin{equation}\label{eqp}
p := \sup_{x \in (0, 1]} \frac{h_a'(x)x}{h_a(x)}, \quad
\Phi_p(x) := \frac{h_a(x)}{x^p}, \quad x > 0.
\end{equation}

\begin{lemma}\label{lemp}
\begin{enumerate}
\item $0 < p \leq 2e^{-1} < 1$.
\item The function $\Phi_p(x)$ is nonincreasing on $(0, 1)$.
\end{enumerate}

\end{lemma}
\begin{proof}
1. Let $\mu := 1/(4a)$. Then\\
$
2\sqrt{\pi}h_a(x)  =  e^{-\mu x^2} \left( 1 + e^{\mu x^2 - \mu x^{-2}} \right) > e^{-\mu x^2},
$\\
$2\sqrt{\pi}h_a'(x)x
=  e^{-\mu x^2}2\mu\left(- x^2 + x^{-2} e^{\mu x^2 - \mu x^{-2}}\right)  \leq e^{-\mu x^2}2 \mu(x^{-2}-x^2)e^{-\mu(x^{-2}-x^2)}
\leq 2e^{-1} e^{-\mu x^2}.
$
Hence, $h_a'(x)x < 2e^{-1} h_a(x)$ and $p \leq 2e^{-1}$. The inequality $p > 0$ is obvious.
\medskip

2. Relations \eqref{eqp} imply the following for $x \in (0, 1]$:
$$
\Phi'_p(x) = \frac{h'_a(x)x^p-h_a(x)px^{p-1}}{x^{2p}} = \frac{1}{x^{p+1}}(h'_a(x)x - ph_a(x)) \leq 0.
$$
Hence, $\Phi_p$ is nonincreasing on $(0, 1)$.
\end{proof}

For any $N \in \bbN$, set
\begin{equation}\label{eqDp12}
D_{p1}  := 2^{\frac{1-p}{2}}-1, \quad D_{p2}  := N\lr{\lr{1+\frac{1}{N}}^{\frac{1+p}{2}}-1}.
\end{equation}

\begin{lemma}\label{lemDp}
For any $N \in \bbN$ we have for $k \geq 1$, $n \geq N$ the following inequalities hold
\begin{equation}\label{ineqDp12}
\left(1+\frac{1}{k}\right)^{\frac{1-p}{2}} \geq 1+ D_{p1}\frac{1}{k},
\quad
\left(1+\frac{1}{n}\right)^{\frac{1+p}{2}} \geq 1+ D_{p2}\frac{1}{n},
\end{equation}
\begin{equation}\label{eqAdd6.1}
D_{p1} < D_{p2}, \quad D_{p1} < \frac{1}{2}.
\end{equation}

\end{lemma}
\begin{proof}
Fix $\al \in (0, 1)$. Note that
\begin{equation}\label{S521}
(1+x)^{1-\al} \leq 1 + x(1-\al), \quad x \in (0, 1].
\end{equation}
Consider the function $P_{\al}(x) := ((1+x)^{\al}-1)/x$. Inequality \eqref{S521} implies that
$$
P'_{\al}(x) = \frac{\al(1+x)^{\al-1}x - ((1+x)^{\al}-1)}{x^2} \leq 0.
$$
Hence, $P_{\al}(x)$ is nonincreasing and, for any $k_1 \geq k_2$, the following holds:
\begin{equation}\label{S522}
P_{\al}(1/k_1) \geq P_{\al}(1/k_2).
\end{equation}
Inequalities \eqref{ineqDp12}  are straightforward consequences of \eqref{S522}
for $\al = (1-p)/2$ and  $\al = (1+p)/2$, respectively. Obviously, $D_{p2} \geq 2^{(1+p)/2} - 1 > D_{p1}$ and $D_{p1}< \sqrt{2}-1<1/2$.
\end{proof}

Set ($D_{p1}$, $D_{p2}$ are given by \eqref{eqDp12})
\begin{equation}
\varkappa  := D_{p2} - D_{p1} - 2(D_{p1}+ D_{p2})D_{p1}\frac{1}{N}. \label{eqep}
\end{equation}
Lemma \ref{lemp} and relations \eqref{eqDp12} imply that $\varkappa > 0$ for large enough $N$.
In what follows, we fix
\begin{equation}\label{eqkap0}
N \in \bbN \quad \mbox{such that} \quad \varkappa > 0.
\end{equation}
Assume that $n_0 = n_0(E)$ satisfies the following. (We remind that the complete list of requirements determining $n_0(E)$ is given in Section~\ref{refAppendixSubsectionAssumptionsN0}.)
\begin{requirement}\label{req1}
$n_0 \geq N$.
\end{requirement}

Set
\begin{equation}\label{eqCn}
C_n := -cn^2 +(h_1 - 2c)an^2 - h_1 \sum_{|k| \leq n-1} y_{n-k}(a(n^2-k^2)).
\end{equation}
In other words, the values $C_n$ are obtained by formally substituting $t = an^2$, $t_k = ak^2$, and $S(t) = \{-(n-1), \dots, n-1\}$ in \eqref{eqSolu}.
In Section \ref{secCn} below, we will prove the following.
\begin{prp}\label{lemCn} There exist $K, K'>0$ such that for $n \geq N$ the following inequalities hold:
\begin{equation}
\left|C_{n+1} - C_n\right|  \leq K\frac{1}{\sqrt{n}}, \quad |C_{n}|  \leq K'\sqrt{n}. \label{eqCn1}\\
\end{equation}
\end{prp}
Fix $K$ and $K'$ from Proposition \ref{lemCn} and set
\begin{equation}
E_0  := \frac{K + \varkappa K'}{(h_1 - 2c)\varkappa}. \label{eqE0}
\end{equation}
Note that $E_0 > 0$ due to \eqref{eqkap0}. For each $E> E_0$, set
\begin{equation} \label{eqaminmax}
\amin_n := a - \frac{2E\sqrt{n}}{2n-1}, \quad \amax_n := a + \frac{2E\sqrt{n}}{2n-1}.
\end{equation}
We assume that $n_0 = n_0(E)$ satisfies the following requirement.
\begin{requirement}\label{ass1.1}
$\amin_n \geq \tau_0/(2n-1)$ for all $n \geq n_0$, where $\tau_0$ was fixed in \eqref{eqtau0}.
\end{requirement}
Note that Requirement \ref{ass1.1} implies that
\begin{equation}\label{eqtau01}
a(n^2-k^2) - 2E\sqrt{n} \geq \tau_0, \quad n \geq n_0, \; |k| \leq n-1.
\end{equation}

Set
\begin{equation} \label{eqdeltan}
\delta_n := \frac{2 A_1}{aD_a}\frac{1}{n} + \frac{ 2A_2^*E
a^{1/2}}{D_a (\amin_n)^{3/2}}
\frac{1}{\sqrt{n}}.
\end{equation}
Obviously,
$$
\amin_n \to a, \quad \delta_{n} \to 0 \quad \mbox{as $n \to \infty$},
$$
$$
\amin_{n+1} > \amin_n, \quad \delta_{n+1} < \delta_{n}, \quad n \geq n_0.
$$
We assume that $n_0 = n_0(E)$ satisfies the following  requirement.
\begin{requirement}\label{ass1.2}
$\delta_n + \delta_{n+1} \leq 2D_{p1}$ for all $n \geq n_0$.
\end{requirement}

Below we will use the following constants $S_{\al}$, $T_{\al}$, and $R$.
For $\al > 0$, let $S_{\al}$ be the smallest number satisfying the inequalities
\begin{equation}\label{eqSal}
\sum_{|k|\leq n-1}\frac{(n-k)^{\al}}{(n^2-k^2)^{\al+1/2}} \leq S_{\al}\frac{1}{\sqrt{n}}, \quad n \geq N.
\end{equation}
For $\al > 1$, let $T_{\al}$ be the smallest number satisfying the inequalities
\begin{equation}\label{eqTalpha}
\sum_{|k|\leq n-1}\frac{1}{(n^2-k^2)^{\al}} \leq T_{\al}\frac{1}{n^{\al}},\quad n \geq N.
\end{equation}
Let $R$ be the smallest number satisfying the inequalities
\begin{equation}\label{eqR}
\sum_{|k|\leq n-1}\frac{1}{(n^2-k^2)^{1/2}} \leq R \; (= \pi), \quad n \geq N.
\end{equation}

\subsection{Candidates for switching moments $t_n$}

\subsubsection{Formulation of a theorem on existence of the candidates}

In this section, we will prove the following result.
\begin{theorem}\label{thqn}
Let the assumptions of Theorem \str{\ref{thMainResult}} hold. Then there exists
a sequence $t_k$, $k\in\bbZ$, such that $t_0,t_{\pm1},\dots,t_{\pm n_0}$ are given by~\eqref{eqEn+ScondtkForuntn},
$$
t_n=an^2+q_n,\quad |q_n| \leq E\sqrt{n},\quad t_{-n}=t_n \quad
\text{for } n \ge n_0+1,
$$
$$
t_{k}<t_{n_0}<t_{n_0+1}<\dots, \quad k = 0, 1, \dots, n_0 - 1,
$$
and the functions
\begin{equation}\label{eqvInsteadOfu}
v_n(t):=-cn^2+(h_1-2c)t-h_1\sum\limits_{|k|\le n-1}y_{n-k}(t-t_k),
\quad n\ge n_0+1,
\end{equation}
satisfy
\begin{equation}\label{equntn0}
v_n(t_n)=0, \quad n \geq n_0 + 1.
\end{equation}
\end{theorem}

\begin{remark}\label{remtnRealSwitch}
The sequence $t_k$ in Theorem \ref{thqn} is a sequence of candidates
for \textit{switching} moments in the following sense. Assume that,
for some $n\ge n_0+1$, we know the following (this is what we will
in particular prove in Sections~\ref{secNablaUn} and~\ref{secEstimatesUn1}
below):
\begin{enumerate}
\item  the nodes $u_0, \dots, u_{n-1}$ switch at time moments
$t_0,\dots,t_{n-1}$, respectively,

\item the nodes $u_n, u_{n+1},u_{n+2},\dots$ do not switch on the
time interval $[0,t_n)$.
\end{enumerate}
Then $v_n(t)$ coincides with the solution $u_n(t)$ of
problem~\eqref{equ_nEquation1}--\eqref{nonideal} on the time
interval $[0,t_n)$ and equality~\eqref{equntn0} implies that $t_n$
is the switching moment of $u_n(t)$.
\end{remark}

\subsubsection{Proof of Theorem \ref{thqn}}

First, we substitute \eqref{eqvInsteadOfu} into \eqref{equntn0}, replace $t_k$ and $t_n$ by $ak^2+ q_k$ and $an^2+ q_n$, respectively, and expand $y_{n-k}$ into the Taylor series around $a(n^2-k^2)$. This yields
\begin{align}
0 = & -cn^2 + (h_1 - 2c)(an^2+q_n) - h_1 \sum_{|k|\leq n-1} y_{n-k}(a(n^2-k^2)) \notag \\
& - h_1 \sum_{|k| \leq n-1} \dot{y}_{n-k}(a(n^2-k^2))(q_n - q_k) \notag \\
& - h_1 \sum_{|k| \leq n-1} \frac{1}{2}\ddot{y}_{n-k}(a(n^2-k^2) + \xi_{n, k})(q_n - q_k)^2 , \quad n \geq n_0 + 1, \label{PAdd1s1}
\end{align}
where $|\xi_{n, k}| \leq |q_n - q_k|$. We introduce the notation
%
\begin{equation}\label{P1s3}
\al_{n, k} := \dot{y}_{n-k}(a(n^2-k^2)), \quad \beta_{n, k}(q_n) := \frac{1}{2}\ddot{y}_{n-k}(a(n^2-k^2)+ \xi_{n, k})(q_n-q_k),
\end{equation}
where we omit an explicit indication of the dependence of $\be_{n, k}$ on $q_k$ with $|k| \leq n-1$. Further, set for $k = 1, 2, \dots, n-1$
\begingroup
\allowdisplaybreaks
\begin{align}
J_{n, k}(q_n) & := \al_{n, k} + \al_{n, -k} + \beta_{n, k}(q_n) + \beta_{n, -k}(q_n), \quad J_{n, 0}(q_n)  := \al_{n, 0} + \beta_{n, 0}(q_n), \label{PAdd1ss3} \\
J_n(q_n) & := \sum_{k = 0}^{n-1} J_{n, k}, \quad D_n(q_n)  := h_1 - 2c -h_1J_n(q_n). \label{PAdd1ss12}
\end{align}
\endgroup
Using this notation and recalling the definition of the constants $C_n$ in \eqref{eqCn}, we rewrite \eqref{PAdd1s1} as follows
(it will also be convenient to replace $n$ by $n+1$):
\begin{equation}\label{PAdd2s2}
C_{n+1} + D_{n+1}(q_{n+1})q_{n+1} + h_1 \sum_{k = 0}^n J_{n+1, k}(q_{n+1})q_k = 0, \quad n \geq n_0.
\end{equation}
Thus, it remains to find a sequence $q_k$, $k \in \mathbb{Z}$, such that $|q_k|\leq E\sqrt{n_0}$ for $k = 0, \pm 1, \dots, \pm n_0$, $|q_{n+1}| \leq E \sqrt{n+1}$, $q_{-(n+1)} = q_{n+1}$ for $n = n_0, n_0 +1, \dots$, and the equalities \eqref{PAdd2s2} hold.

First, we note that $q_0, \dots, q_{\pm n_0}$ are already prescribed by the assumption of the theorem. Moreover, \eqref{PAdd2s2} holds with $n+1$ replaced by $n_0$:
\begin{equation}\label{PAdd2s3}
C_{n_0} + D_{n_0}(q_{n_0})q_{n_0} + h_1 \sum_{k = 0}^{n_0-1}J_{n_0, k}(q_{n_0})q_k = 0.
\end{equation}
Indeed, Requirement \ref{ass1.1} implies that $t_{n_0} > t_k$, $k = 0, \dots, n_0 - 1$. Therefore, for all $t \in [t_{n_0-1}, t_{n_0})$ holds
$
S(t) = \{ -(n_0 - 1), \dots, n_0 - 1 \},
$
in \eqref{eqSolu}  and
\begin{equation}\label{PAdd2s4}
u_{n_0}(t_{n_0}) = -c n_0^2 + (h_1 - 2c)t_{n_0} - h_1 \sum_{|k| \leq n_0 - 1} y_{n_0 - k}(t_{n_0} - t_k).
\end{equation}
Hence, \eqref{PAdd2s3} is obtained in the same way as \eqref{PAdd2s2} from \eqref{eqvInsteadOfu} and \eqref{equntn0}.

Now we proceed by induction. Fix $n \geq n_0$. Suppose, we have constructed the desired sequence $q_0, \dots, q_n$. Let us find $q_{n+1}$ satisfying $|q_{n+1}| \leq E \sqrt{n+1}$ and equation \eqref{PAdd2s2}.
We rewrite equation~\eqref{PAdd2s2} in the form
\begin{equation}\label{eq2FixedPoint}
q_{n+1} = \bF(q_{n+1}),\qquad \bF(q_{n+1}):=-
\dfrac{C_{n+1}}{D_{n+1}(q_{n+1})}-h_1\sum\limits_{k=0}^n \dfrac{J_{{n+1}, k}(q_{n+1})
}{D_{n+1}(q_{n+1})} \cdot q_k.
\end{equation}

To prove Theorem~\ref{thqn}, it now suffices to show that if
$q_k\in[-E\sqrt{n},E\sqrt{n}]$ for $k=0,\pm 1,\dots,\pm n$, then  $\bF$
has a fixed point on the interval $[-E\sqrt{n+1},E\sqrt{n+1}]$.

To do so, we need to show that $\bF$ maps the interval
$[-E\sqrt{n+1},-E\sqrt{n+1}]$ into itself. Let us indicate the
main difficulty on this way. We will see in Sections \ref{secCn} and~\ref{secJ}
that $C_n \sim \sqrt{n}$, $D_n \sim 1/\sqrt{n}$, and
$J_{n+1,k}(q_{n+1})\sim 1/\sqrt{n^2-k^2}$, provided that $|q_{n+1}| \leq E\sqrt{n+1}$.
Therefore, the straightforward attempt to estimate $|\bF(q_{n+1})|$ would yield
\begin{equation}\label{ip12s}
|\bF(q_{n+1})|\le
\left|\dfrac{C_{n+1}}{D_{n+1}(q_{n+1})}\right|+E\sqrt{n}\cdot
h_1\sum\limits_{k=0}^n \left|\dfrac{J_{{n+1},
k}(q_{n+1})}{D_{n+1}(q_{n+1})}\right|,
\end{equation}
and we would obtain nothing better than $|\bF(q_{n+1})|\le
\const\cdot n$.

To overcome this difficulty, we will use the following trick. Note
that, by the induction hypothesis, \eqref{PAdd2s2} holds with $n+1$
replaced by $n$. Therefore, we can multiply \eqref{PAdd2s2} by
$1+\varkappa/n$ with an appropriate $\varkappa>0$ and subtract
\eqref{PAdd2s2} with $n+1$ replaced by $n$. As a result, we will obtain
the equation
\begin{equation}\label{ip12ss}
q_{n+1} = \tilde{\bF}(q_{n+1}),\qquad \tilde{\bF}(q_{n+1})= -\dfrac{\tilde C_{n+1}}{\tilde D_{n+1}(q_{n+1})} -
h_1\sum\limits_{k=0}^n  \dfrac{\tilde J_{{n+1},
k}(q_{n+1})}{\tilde D_{n+1}(q_{n+1})}\cdot q_k,
\end{equation}
which is equivalent to \eqref{eq2FixedPoint}. The advantage of this
new representation will be that we will obtain $\tilde C_n \sim
1/\sqrt{n}$ and  $\tilde D_n \sim 1/\sqrt{n}$. Hence, the first term
in the formula for $\tilde{\bF}$ can be estimated by a
constant $\alpha_1>0$. Furthermore, we will show that the expression
$h_1\sum\limits_{k=0}^n \left|\dfrac{\tilde J_{{n+1},
k}(q_{n+1})}{\tilde D_{n+1}(q_{n+1})}\right|$ is estimated by
$1-\alpha_2/\sqrt{n}$ with $\alpha_2>0$. Therefore, \eqref{ip12ss}
will yield
\begin{equation}\label{ip12sss}
|\tilde{\bF}(q_{n+1})|\le \alpha_1+E\sqrt{n}-E\alpha_2<
E\sqrt{n+1},
\end{equation}
if $E\ge E_0=\alpha_1/\alpha_2$. In particular, it will
turn out that the appropriate $\varkappa$ is given by \eqref{eqep}
and $E_0$ by \eqref{eqE0}. Interestingly, $\varkappa=0$ would not be
sufficient for this scheme as it would  then  follow that
$\alpha_2=0$.

To make the above argument rigorous, we need the following proposition, in which we do not explicitly
indicate the dependence of the functions on $q_{n+1}$.

\begin{prp}\label{lEstimatesJFixedPoint}
Let the assumptions of Theorem \str{\ref{thMainResult}} hold. Then for any $n \geq n_0$ and  $q_0, q_{\pm 1}, \dots, q_{\pm (n+1)} \in [-E\sqrt{n+1},E\sqrt{n+1}]$, the following holds with $\varkappa$ given by \eqref{eqep}$:$
\begin{enumerate}
\item $J_{n, k} \geq 0$, $k = 0, 1, \dots, n-1$,

\item $J_{n, k} - \left(1 + \dfrac{\varkappa}{n}\right)J_{n+1, k}  \geq 0$, $k = 0, 1, \dots, n-1$,

\item $D_n - h_1 \left(1 + \dfrac{\varkappa}{n}\right) J_{n+1, n}  \geq 0$.

\end{enumerate}
\end{prp}
Now, assuming that Proposition~\ref{lEstimatesJFixedPoint} is true, we complete the proof of Theorem~\ref{thqn}. After that, in Section \ref{secJ}, we prove Proposition~\ref{lEstimatesJFixedPoint}.

Substituting $n$ instead of $n+1$  into \eqref{PAdd2s2} and using \eqref{PAdd2s3} for $n = n_0$ or the induction hypothesis for $n > n_0$ (and omitting the dependence of $J_{n, k}$ on $q_n$), we obtain the equation
$$
C_n + D_n q_n + h_1 \sum_{k = 0}^{n-1}J_{n, k}q_k = 0.
$$

Multiplying \eqref{PAdd2s2} by $1+ \varkappa/n$ and subtracting the latter expression we have
\begin{align*}
0 & = (C_{n+1}-C_n) + \frac{\varkappa}{n}C_{n+1} + \lr{1+ \frac{\varkappa}{n}}D_{n+1}q_{n+1} \\
& \;\;\;\; +  \lr{\lr{1 +\frac{\varkappa}{n}}h_1J_{n+1, n} - D_n } q_n + h_1\sum_{k = 0}^{n-1}\lr{\lr{1 +\frac{\varkappa}{n}}J_{n+1, k} - J_{n, k}}q_k.
\end{align*}
Equivalently (cf. \eqref{ip12ss}), $q_{n+1} = \tilde{\bF}(q_{n+1})$, where $\tilde{\bF}(q_{n+1})$ satisfies
\begin{align*}
 \lr{1+ \frac{\varkappa}{n}}D_{n+1}\tilde{\bF}(q_{n+1}) & =  -\lr{(C_{n+1}-C_n) + \frac{\varkappa}{n}C_{n+1}}  + \lr{D_n - \lr{1 +\frac{\varkappa}{n}}h_1 J_{n+1, n}} q_n  \\
  & \;\;\;\; + h_1\sum_{k = 0}^{n-1}\lr{J_{n, k} - \lr{1 +\frac{\varkappa}{n}}J_{n+1, k}}q_k.
\end{align*}
According to Proposition \ref{lEstimatesJFixedPoint}, all the coefficients at $q_k$, $q_n$, $q_{n+1}$ are positive. The inductive hypothesis $|q_k| \leq E\sqrt{n}$, $k = 0, \dots, n$, Proposition \ref{lemCn}, and the inequality $E \geq E_0$ imply that
\begin{multline*}
 \left(1+ \frac{\varkappa}{n}\right)D_{n+1}|\tilde{\bF}(q_{n+1})|  \leq \left|(C_{n+1}-C_n) + \frac{\varkappa}{n}C_{n+1}\right| \\
  + \left(D_n - \left(1\! +\!\frac{\varkappa}{n}\right)h_1J_{n+1, n}\right) E\sqrt{n} +
   h_1\sum_{k = 0}^{n-1}\left(J_{n, k} - \left(1 +\frac{\varkappa}{n}\right)J_{n+1, k}\right)E\sqrt{n}  \\
     \leq \frac{K + \varkappa K'}{\sqrt{n}} + E\sqrt{n}\left(\!\!\left(\!D_n + h_1\! \sum_{k = 0}^{n-1} J_{n, k}\!\right)\!\! -\!h_1\!\!\left(1\!+\! \frac{\varkappa}{n}\right)\!\!\left(\!J_{n+1, n}\!+\! \sum_{k = 0}^{n-1}J_{n+1, k}\!\right)\!\!\right)\!\!.
\end{multline*}
Combining this with \eqref{PAdd1ss12}  yields (cf. \eqref{ip12sss})
\begin{align*}
 \left(1+ \frac{\varkappa}{n}\right)D_{n+1}|\tilde{\bF}(q_{n+1})|  & \leq  (K + \varkappa K')\frac{1}{\sqrt{n}} + E\sqrt{n}\left((h_1-2c) - \lr{1+ \frac{\varkappa}{n}}((h_1-2c) - D_{n+1})\right) \\
   & =  (K + \varkappa K')\frac{1}{\sqrt{n}} + E\sqrt{n}D_{n+1}\left(1+ \frac{\varkappa}{n}\right) - E\sqrt{n}(h_1-2c)\frac{\varkappa}{n}.
\end{align*}
The latter estimate and the inequality $E \geq E_0$, where $E_0$ is given by \eqref{eqE0}, imply
\begin{equation*}
\left(1+ \frac{\varkappa}{n}\right)D_{n+1}|\tilde{\bF}(q_{n+1})|  \leq  E\sqrt{n}D_{n+1}\lr{1+ \frac{\varkappa}{n}}.
\end{equation*}
Hence, $|\tilde{\bF}(q_{n+1})| \leq E\sqrt{n} < E\sqrt{n+1}$.
Therefore, $\tilde{\bF}$ maps $[-E\sqrt{n+1}, E\sqrt{n+1}]$ into itself.

Furthermore, the function $\tilde{\bF}$ is continuous, because the
functions
$$
\beta_{m, k}(q_m) = \begin{cases}
\frac{y_{m-k}(a(m^2-k^2) + q_m-q_k) -
y_{m-k}(a(m^2-k^2))}{q_m-q_k} - \al_{m, k}, \; & q_m \ne q_k,\\
0, \; & q_m = q_k,
\end{cases}
$$
are continuous with respect to  $q_m \in \mathbb{R}$ for $k = 0,
\dots, m-1$ and $D_{n+1}(q_{n+1}) > 0$ by
Proposition~\ref{lEstimatesJFixedPoint}. Hence, by Brouwer's
fixed-point theorem, the map $\tilde{\bF}$ has a fixed point on the interval
$[-E\sqrt{n+1}, E\sqrt{n+1}]$. The latter implies Theorem~\ref{thqn}.

It remains to prove Propositions \ref{lemCn} and
\ref{lEstimatesJFixedPoint}, which we do in the next two sections.

\subsection{Proof of Proposition~$\ref{lemCn}$}\label{secCn}

Taking into account equalities \sref{eqr01t} and \eqref{eqah1f} we write \eqref{eqCn} as follows:
\begingroup
\allowdisplaybreaks
$$
C_n = C_n^{(1)} + C_n^{(2)} + C_n^{(3)},
$$
where
\begin{align*}
C_n^{(1)} & := h_1\lr{n^2 I_F - \sum_{|k| \leq n-1} \sqrt{a(n^2-k^2)}f\lr{\frac{n-k}{\sqrt{a(n^2-k^2)}}}}, \label{Add24.1}\\
C_n^{(2)} & := -h_1 \sum_{|k| \leq n-1}\frac{1}{\sqrt{a(n^2-k^2)}}\tilde{f}\lr{\frac{n-k}{\sqrt{a(n^2-k^2)}}},\notag\\
C_n^{(3)} & := -h_1 \sum_{|k| \leq n-1} \tilde{r}_1(n-k, a(n^2-k^2)).\notag
\end{align*}
\endgroup

\begin{proof}[Proof of the second inequality in $\sref{eqCn1}$]
Below we separately estimate $C_{n}^{(1)}$ and $C_{n}^{(2)} + C_{n}^{(3)}$.

\noindent\textbf{Step 1.}
Proposition \ref{lemIntSqrt1} (item 1) applied to the function $F(x)$ (see \eqref{eqFG}), implies that
\begin{equation}\label{eqK1p}
|C_n^{(1)}| \leq h_1 L_1(F, N) \sqrt{n} = K'_1\sqrt{n}, \quad K'_1 := h_1 L_1(F, N).
\end{equation}

\noindent\textbf{Step 2.} Due to \eqref{eqr01t}, \eqref{eqr01t}, and \eqref{eqA},
$$
|C_n^{(2)} + C_n^{(3)}| = \left|h_1 \!\! \sum_{|k|\leq n-1}r_0(n-k, a(n^2-k^2))\right| \leq h_1  \!\!\sum_{|k| \leq n-1} \frac{A_0}{\sqrt{a(n^2-k^2)}} \leq  l_1,
$$
where, by \eqref{eqR},
\begin{equation}\label{eql1}
l_1 := h_1 A_0 \frac{1}{a^{1/2}}R.
\end{equation}
Hence inequality \sref{eqCn1} is satisfied for
\begin{equation}\label{eqK'}
K' := K'_1+ l_1\frac{1}{\sqrt{N}}.
\end{equation}
\end{proof}

\begin{proof}[Proof of the first inequality in $\sref{eqCn1}$]
Below we separately estimate $C_{n+1}^{(1)} - C_n^{(1)}$, $C_{n+1}^{(2)} - C_n^{(2)}$, and $C_{n+1}^{(3)} - C_n^{(3)}$ in Steps 1, 2, and 3 respectively.

\noindent\textbf{Step 1.}
Applying Proposition \ref{lemIntSqrt1} (item 2) to the function $F(x)$, we conclude that
\begin{equation}\label{P16s}
|C_{n+1}^{(1)} - C_n^{(1)}| \leq K_1 \frac{1}{\sqrt{n}}, \quad K_1 := h_1 \bar{L}_1(F, N).
\end{equation}

\noindent\textbf{Step 2.} Set
\begin{equation}\label{eqtF}
\tilde{F}(x) := \frac{1}{\sqrt{a(1-x^2)}}\tilde{f}\lr{\frac{1-x}{\sqrt{a(1-x^2)}}}, \quad I_{\tilde{F}} := \int_{-1}^1 \tilde{F}(x) dx.
\end{equation}
Note that $\lim \limits_{x \to 1} \tilde{F}(x)(1-x)^{1/2} < 0$. Therefore, by Proposition \ref{lemIntSqrt2} applied to the function $-\tilde{F}(x)$, for some constants $L_{\tilde{f}}$, $L_{\tilde{f}}^*$, $l_{\tilde{f}} \geq 0$, we have
\begin{equation}\label{eqLt}
L_{\tilde{f}}^* \frac{1}{\sqrt{n}} - l_{\tilde{f}}\frac{1}{n} \leq -I_{\tilde{F}} - \frac{1}{h_1}C_n^{(2)} \leq L_{\tilde{f}} \frac{1}{\sqrt{n}}.
\end{equation}
Hence,
\begin{equation}\label{P16ss}
|C_{n+1}^{(2)} - C_n^{(2)}| \leq K_2\frac{1}{\sqrt{n}}, \quad K_2 := h_1 \lr{ L_{\tilde{f}} - L_{\tilde{f}}^* + l_{\tilde{f}}\frac{1}{\sqrt{N}}}
\end{equation}

\noindent\textbf{Step 3.}
Inequalities \sref{eqA} and \eqref{eqTalpha} imply that
$$
|C_n^{(3)}| \leq h_1\tilde{A}_1 \sum_{|k| < n}\frac{1}{(a(n^2-k^2))^{3/2}} \leq h_1 \tilde{A}_1 \frac{1}{a^{3/2}}T_{3/2}\frac{1}{n^{3/2}}.
$$
Hence,
\begin{equation}\label{P16sss}
|C_{n+1}^{(3)} - C_n^{(3)}| \leq K_3\frac{1}{n^{3/2}}, \quad K_3 := 2h_1\tilde{A}_1 \frac{1}{a^{3/2}}T_{3/2}
\end{equation}

Summarising \eqref{P16s}, \eqref{P16ss}, and \eqref{P16sss} yields
\begin{equation}\label{eqK}
|C_{n+1} - C_n| \leq K \frac{1}{\sqrt{n}}, \quad K := K_1 + K_2 + K_3 \frac{1}{N}.
\end{equation}

\end{proof}

\subsection{Proof of Proposition~$\ref{lEstimatesJFixedPoint}$}\label{secJ}

\subsubsection{Proof of Proposition~\erm{\ref{lEstimatesJFixedPoint}}\erm{:} Preliminaries}\label{secJprelim}

For $k = 1, \dots, n-1$, consider the following representation of $J_{n, k}$ (see \eqref{PAdd1ss3} and \eqref{eqr12}):
$$
J_{n, k} = J_{n, k}^{\rm main} + (w_{n, k} + w_{n, -k}) + (\be_{n, k}+
\be_{n, -k}),
$$
$$
\Jnkm := \gamma_{n, k} + \gamma_{n, -k}, \quad \gamma_{n, \pm k} := \frac{1}{\sqrt{a(n^2-k^2)}}h\left(\sqrt{\frac{n \mp k}{a(n \pm k)}}\right), \quad
w_{n, \pm k} := r_1(n \mp k, (a(n^2-k^2)).
$$
For $k=0$, we have
$$
J_{n, 0} = J_{n, 0}^{\rm main} + w_{n, 0} + \be_{n, 0}, \quad
J_{n, 0}^{\rm main} := \gamma_{n, 0} := \frac{1}{\sqrt{an^2}} h\left(\frac{1}{\sqrt{a}}\right), \quad w_{n, 0} := r_1(n, an^2).
$$

The general idea is to prove each assertion of Proposition \ref{lEstimatesJFixedPoint} for $\Jnkm$ first, and then consider $\Jnk$ as a small perturbation of $\Jnkm$. We formulate this fact as a lemma.

\begin{lemma}\label{lemdeltan} Let $\delta_n$ be given by \eqref{eqdeltan}. Then
\begin{equation}\label{eqDeltaJnk}
\frac{|J_{n, k}-J_{n, k}^{\rm main}|}{J_{n, k}^{\rm main}} \leq \delta_n \frac{n}{n+k} \frac{1}{n-k}, \quad k \in \{0, 1, \dots, n-1\}.
\end{equation}
\end{lemma}

\begin{proof}
Fix $k \in \{1, \dots, n-1\}$ (the case $k = 0$ is similar).
Due to \eqref{eqA} and \eqref{eqtau01},
\begin{equation}\label{P17s1}
|w_{n, k} + w_{n, -k}| \leq 2 A_1 \frac{1}{(a(n^2-k^2))^{3/2}}.
\end{equation}
Definition \eqref{P1s3} of $\be_{n, k}$, \eqref{eqtau01}, \sref{eqAB2*} and \eqref{eqaminmax} imply that
\begin{align}
|\be_{n, k} + \be_{n, -k}| & \leq \max_{|\xi_{n, k}| \leq |q_n - q_k|}|\ddot{y}_{n-k}(a(n^2-k^2) + \xi_{n, k})| \cdot |q_n-q_k| \notag \\
& \leq  \frac{2A_2^*E\sqrt{n}}{(a(n^2-k^2) - 2E\sqrt{n})^{3/2}} \leq \frac{2A_2^*E\sqrt{n}}{(\amin_n(n^2-k^2))^{3/2}} \label{P17s2}.
\end{align}

Note that, by \sref{eqfDa},
\begin{equation}\label{P18s1}
J_{n, k}^{\rm main} \geq \frac{D_a}{(a(n^2-k^2))^{1/2}}.
\end{equation}

Relations \eqref{P17s1}--\eqref{P18s1} imply
\begingroup
\allowdisplaybreaks
\begin{align*}
\frac{|w_{n, k}+ w_{n, -k}|}{J_{n, k}^{\rm main}} & \leq  \frac{2
A_1\frac{1}{(a(n^2-k^2))^{3/2}}}{\frac{D_a}{(a(n^2-k^2))^{1/2}}}
 =   \frac{2 A_1}{aD_a}\frac{1}{n+k} \frac{1}{n-k},\\
\frac{|\be_{n, k} + \be_{n, -k}|}{J_{n, k}^{\rm main}} & \leq \frac{
 \frac{2A_2^*E\sqrt{n}}{(\amin_n(n^2-k^2))^{3/2}}}
{\frac{D_a}{(a(n^2-k^2))^{1/2}}} =   \frac{2
A_2^* E a^{1/2}}{D_a(\amin_n)^{3/2}}
\frac{\sqrt{n}}{n+k} \frac{1}{n-k}.
\end{align*}
\endgroup
Hence, inequality \sref{eqDeltaJnk} holds.
\end{proof}

\subsubsection{Proof of Proposition~$\ref{lEstimatesJFixedPoint}$\erm{:} Part \erm{1}}

Since $\Jnkm > 0$, it suffices to show that the right-hand side in \eqref{eqDeltaJnk} is less than or equal to~1. The latter is true because $\delta_n \leq 1$ due to \eqref{eqAdd6.1} and Requirement~\ref{ass1.2}.

\subsubsection{Proof of Proposition~$\ref{lEstimatesJFixedPoint}$\erm{:} Part \erm{2}}

Fix $k \in \{ 1, \dots, n-1 \}$ (the case $k = 0$ can be treated similarly).
Note that
\begin{equation*}
\Jnkm  = \frac{1}{\sqrt{a(n^2-k^2)}} \left(\sqrt{\frac{n-k}{n+k}}\right)^p \Phi_p\left(\sqrt{\frac{n-k}{n+k}}\right)
 =  \frac{\Phi_p\left(\sqrt{\frac{n-k}{n+k}}\right)}{\sqrt{a}(n-k)^{(1-p)/2}(n+k)^{(1+p)/2}}
,
\end{equation*}
where $\Phi_p(x)$, $x \in (0, 1]$, is given by \eqref{eqp}. Taking into account that $\Phi_p(x)$ is nonincreasing (see Lemma \ref{lemp}, item 2) and using \eqref{ineqDp12}, we have
\begin{align}
\frac{\Jnkm}{J_{n+1, k}^{\rm main}} & \geq
\frac{(n+1-k)^{(1-p)/2}(n+1+k)^{(1+p)/2}}{(n-k)^{(1-p)/2}(n+k)^{(1+p)/2}} = \left( 1+ \frac{1}{n-k}\right)^{(1-p)/2} \left(1+ \frac{1}{n+k}\right)^{(1+p)/2} \notag \\
& \geq \left(1+ D_{p1}\frac{1}{n-k}\right)\left(1+ D_{p2}\frac{1}{n+k}\right) \geq
1+ D_{p1}\frac{1}{n-k} + D_{p2}\frac{1}{n+k}.\label{eqAdd10.1}
\end{align}

Combining \eqref{eqAdd10.1} with \sref{eqDeltaJnk} yields
\begin{align*}
\frac{J_{n, k}}{J_{n+1, k}} & =
\frac{\Jnkm}{J_{n+1, k}^{\rm main}} \frac{ 1 + \frac{J_{n, k} - \Jnkm}{\Jnkm}}{1 + \frac{J_{n+1, k} - J_{n+1, k}^{\rm main}}{J_{n+1, k}^{\rm main}}} \geq
\frac{\Jnkm}{J_{n+1, k}^{\rm main}}\frac{1 - \frac{\delta_n}{n-k}\frac{n}{n+k}}{1+\frac{\delta_{n+1}}{n+1-k}\frac{n+1}{n+1+k}} \\
& \geq  \frac{J_{n, k}^{\rm main}}{J_{n+1, k}^{\rm main}}
\left( 1 - \frac{\delta_n}{n-k} \frac{n}{n+k} \right) \left( 1+\frac{\delta_{n+1}}{n-k}\frac{n}{n+k} \right)^{-1} \\
& \geq  \left(1+ D_{p1}\frac{1}{n-k} + D_{p2}\frac{1}{n+k}\right)\left(1- (\delta_n + \delta_{n+1})\frac{n}{n+k}\frac{1}{n-k}\right).
\end{align*}
It is easy to show that Requirenment \ref{ass1.2} implies that the minimum of the last expression is achieved for $k = 0$. Hence, using Requirenment \ref{ass1.2} again, we obtain ($\varkappa$ is given by \eqref{eqep})
\begin{align*}
\frac{J_{n, k}}{J_{n+1, k}} \geq  &
\left(1 + (D_{p1} + D_{p2})\frac{1}{n}\right)\left(1- (\delta_n + \delta_{n+1})\frac{1}{n}\right) \\
\geq & \; 1 + \left(D_{p1}+ D_{p2} - (\delta_n + \delta_{n+1}) - (D_{p1}+ D_{p2})(\delta_n + \delta_{n+1})\frac{1}{N}\right)\frac{1}{n}
\geq  1+ \frac{\varkappa}{n}.
\end{align*}

\subsubsection{Proof of Proposition~$\ref{lEstimatesJFixedPoint}$\erm{:} Part \erm{3}}

Using \eqref{PAdd1ss12} and the notation in the beginning of Section \ref{secJprelim}, we write
\begin{equation}\label{eq19.1}
D_n - h_1 J_{n+1, n}\lr{1+ \frac{\varkappa}{n}} = (h_1 - 2c) - h_1 J_n - h_1 J_{n+1, n} \left(1 +
\frac{\varkappa}{n}\right) = \Sigma_{1, n} - (\Sigma_{2, n} + \Sigma_{3, n} + \Sigma_{4, n}),
\end{equation}
where
\begingroup
\allowdisplaybreaks
\begin{align*}
\Sigma_{1, n} & := (h_1 - 2c) - h_1\sum_{k = 0}^{n-1} \Jnkm,  &
\Sigma_{2, n} & := h_1 \gamma_{n+1, n} \left(1 + \frac{\varkappa}{n}\right), \\
\Sigma_{3, n} & := h_1\sum_{k = 0}^{n-1} \left(\Jnk - \Jnkm \right), &
\Sigma_{4, n} & := h_1 \left( \left(J_{n+1, n} - J_{n+1, n}^{\rm main}\right) + \gamma_{n+1, -n} \right)\left(1 + \frac{\varkappa}{n}\right).
\end{align*}
\endgroup

In steps 1--4 below, we will estimate $\Sigma_{1, n}, \dots, \Sigma_{4, n}$. We will see that $\Sigma_{1, n}$ and $\Sigma_{2, n}$ are ``large'' with respect to $\Sigma_{3, n}$ and $\Sigma_{4, n}$, which motivates the splitting in \eqref{eq19.1}.
Set
$$
\tilde{H}_1(x) := H(x) + H(-x), \quad x \in [0, 1),
$$
where $H(x)$ is given by \eqref{eqHH1}. Then for $k = 1, \dots, n-1$ we have
\begin{equation*}
J_{n, k}^{\rm main}  = \frac{1}{n} \tilde{H}_1\lr{\frac{k}{n}}, \quad J_{n, 0}^{\rm main}  = \frac{1}{2n} \tilde{H}_1(0).
\end{equation*}

\noindent\textbf{Step 1.}
Proposition \ref{lemIntSqrt2} applied to the function $H(x)$ (see \eqref{eqHH1}) and equation \eqref{eqIntfgh} imply
\begin{equation}\notag
\Sigma_{1, n} \geq h_1\lr{C_H \frac{1}{\sqrt{n}} - l_H\frac{1}{n}},
\end{equation}
\begin{equation}\label{eqCH}
C_H := L_2^*(H, N), \quad  l_H := l_2(H, N).
\end{equation}

\noindent\textbf{Step 2.}
Since the function $h(x)$ is decreasing,
$$
\Sigma_{2, n} \leq h_1 \frac{1}{\sqrt{2a}}h(0) \frac{1}{\sqrt{n}}\left(1 + \frac{\varkappa}{n}\right).
$$

\noindent\textbf{Step 3.}
Using \eqref{eqDeltaJnk} and  Proposition \ref{propL3} applied to the function $\dfrac{1}{1-x^2}H(x)$, we have
$$
|\Sigma_{3, n}| \leq h_1 \delta_n \sum_{k = 0}^{n-1}\frac{n}{n^2-k^2}\Jnkm \leq h_1 \delta_n C_{H2}\frac{1}{\sqrt{n}},
$$
where
\begin{equation}\label{eqCH2}
C_{H2} :=  L_3\lr{\frac{1}{1-x^2}H}.
\end{equation}

\noindent\textbf{Step 4.}
Using \eqref{eqDeltaJnk}, we obtain
$$
|J_{n+1, n} - J_{n+1, n}^{\rm main}| \leq  \; \delta_{n+1}\frac{n+1}{2n+1}J_{n+1, n}^{\rm main}
\leq  \; \delta_{n+1}\frac{n+1}{2n+1} \left(\frac{h(0)}{\sqrt{2a}}\frac{1}{\sqrt{n}} + \gamma_{n+1, -n}\right).
$$
Note that
$$
\gamma_{n+1, -n} \leq \frac{1}{\sqrt{2a}} h\left(\sqrt{\frac{2n+1}{a}}\right)\frac{1}{\sqrt{n}}.
$$

Summarising the last two inequalities, we have
$$
\Sigma_{4, n}  \leq h_1 \frac{1}{\sqrt{2a}}\left(1 + \frac{\varkappa}{n}\right)  \lr{h(0) \frac{n+1}{2n+1}\delta_{n+1} +
(1+ \delta_{n+1})h\left(\sqrt{\frac{2n+1}{a}}\right)}\frac{1}{\sqrt{n}}.
$$

Steps 1--4 yield Proposition~$\ref{lEstimatesJFixedPoint}$ (part 2), if $n_0 = n_0(E)$ satisfies the following.
\begin{requirement}\label{ass2}
For $n \geq n_0$, the following holds\str{:}
\begin{multline*}
C_H - \frac{1}{\sqrt{2a}}h(0) \geq l_H\frac{1}{\sqrt{n}} + \frac{\varkappa}{\sqrt{2a}}h(0)\frac{1}{n}  + C_{H2}\delta_n  \\
+ \frac{1}{\sqrt{2a}}\left(1 + \frac{\varkappa}{n}\right)
\left(
h(0)\frac{n+1}{2n+1}\delta_{n+1} +
(1+\delta_{n+1})h\left(\sqrt{\frac{2n+1}{a}}\right)
\right).
\end{multline*}
\end{requirement}
Note that, according to \eqref{eqCH} and Proposition \ref{lemIntSqrt2}, $C_H - h(0)/\sqrt{2a} > 0$.




\section{Asymptotics for $\nabla
u_{n}(t_{n})$}\label{secNablaUn}

We consider the sequence $t_k$ ($k\in\bbZ$) given by
Theorem~\ref{thqn} and the quantities
\begin{equation}\label{eqnablauntn0}
\nabla v_n(t_n)=-2cn-c-h_1\sum\limits_{|k|\le n-1}\nabla
y_{n-k}(t_n-t_k),\quad n\ge n_0+1,
\end{equation}
where $v_n(t)$ is given by~\eqref{eqvInsteadOfu}.

\begin{remark} Under the assumptions of Remark~\ref{remtnRealSwitch}, we have
$
 \nabla v_n(t_n)=\nabla u_n(t_n).
$
\end{remark}

In this section, we will prove that $\nabla v_n(t_n)<0$.

\begin{theorem}\label{thNablaun}
There exists $A_{\nabla}>0$ depending on $h_1,c,E$ such that, for all $n\ge n_0+1$,
$$
\left|\nabla v_n(t_n)+\dfrac{3h_1}{4}\right| \le A_{\nabla} n^{-1/2},\qquad \nabla v_n(t_n) \le -\dfrac{3h_1}{8}.
$$
\end{theorem}

\subsection{Preliminaries}
Set
\begin{equation}\label{eqxnk}
x_{n, k}:=\dfrac{k}{n}.
\end{equation}
Consider a constant $K_{h1}$  such that
\begin{equation}\label{eqConstS11}
\left| \sum\limits_{|k|\le
n-1}\dfrac{1}{n}\dfrac{1}{a(1-x_{n, k}^2)}h'\left(\dfrac{1}{\sqrt{a}}\sqrt{\dfrac{1-x_{n, k}}{1+x_{n, k}}}\right)\right|\le
K_{h1},\quad n\ge N.
\end{equation}
Such a constant exists because the left-hand side
in~\eqref{eqConstS11} is the Riemann sum of a finite integral
(note that $h'(0)=0$).

\subsection{Leading order terms}

Substituting $t_k$ given by Theorem~\ref{thqn}
into~\eqref{eqnablauntn0}, we have
\begin{equation}\label{eqNablauntn1}
\nabla v_n(t_n)=-2cn-c-h_1\sum\limits_{|k|\le n-1}\nabla
y_{n-k}(a(n^2-k^2)+q_n-q_k).
\end{equation}
Due to the Taylor expansion,
\begin{equation}\label{eqNablauntn1'}
\sum\limits_{|k|\le n-1}\nabla
y_{n-k}(a(n^2-k^2)+q_n-q_k)=\sum\limits_{|k|\le n-1}\nabla
y_{n-k}(a(n^2-k^2)) + \Sigma_{1,n},
\end{equation}
where
\begin{equation}\label{eqSigma1}
\Sigma_{1,n}  :=\sum\limits_{|k|\le n-1}\nabla\dot
y_{n-k}(a(n^2-k^2))(q_n-q_k)+\dfrac{1}{2}\sum\limits_{|k|\le n-1}\nabla\ddot
y_{n-k}(a(n^2-k^2)+\xi_{n, k})(q_n-q_k)^2
\end{equation}
with $ |\xi_{n, k}|\le|q_n-q_k|\le 2En^{1/2}.$
Using~\eqref{eqw01} and the functions $G(x)=G(a,x)$ and
$H(x)=H(a,x)$ given by~\eqref{eqFG} and~\eqref{eqHH1}, we represent
the sum in~\eqref{eqNablauntn1'} as follows:
\begin{equation}\label{eqNablauntn2}
\sum\limits_{|k|\le n-1}\nabla y_{n-k}(a(n^2-k^2))
 = n \sum\limits_{|k|\le n-1}\dfrac{1}{n} G(x_{n, k}) +
\dfrac{1}{2}\sum\limits_{|k|\le n-1}\dfrac{1}{n} H(x_{n, k})+
\Sigma_{2,n},
\end{equation}
where
\begin{equation}\label{eqSigma2}
\Sigma_{2,n}:=\sum\limits_{|k|\le n-1} w_0(n-k,a(n^2-k^2)).
\end{equation}

Now we replace the sums in~\eqref{eqNablauntn2} by the
integrals. Set (recall that $G(1) = -1/2$)
\begingroup
\allowdisplaybreaks
\begin{align}
\Sigma_{g,n}&:=\sum\limits_{|k|\le n-1}\dfrac{1}{n}
G(x_{n, k})-I_G+\dfrac{G(1)}{2n}
=\sum\limits_{|k|\le n-1}\dfrac{1}{n} G(x_{n, k})-I_G-\dfrac{1}{4n},\label{eqSigmag}\\
\Sigma_{h,n}&:=\sum\limits_{|k|\le n-1}\dfrac{1}{n}
H(x_{n, k})-I_H,\label{eqSigmah}
\end{align}
\endgroup
where $I_G=I_G(a)$ and $I_H=I_H(a)$ are given by~\eqref{eqIntfgh}. Then~\eqref{eqNablauntn2} takes the form
\begin{equation}\label{eqNablauntn2'}
\sum\limits_{|k|\le n-1}\nabla y_{n-k}(a(n^2-k^2))
 = I_G n +
\dfrac{I_H}{2}+\dfrac{1}{4}+n\Sigma_{g,n}+\dfrac{\Sigma_{h,n}}{2}+\Sigma_{2,n}.
\end{equation}

Combining~\eqref{eqNablauntn1}, \eqref{eqNablauntn1'}, and
\eqref{eqNablauntn2'} and using Lemma~\ref{lema}, we obtain
\begin{align}
\nabla v_n(t_n) &=  (-2c-h_1 I_G)n +
 \left(-c-\dfrac{h_1I_H}{2}-\dfrac{h_1}{4}\right)  -h_1\left(\Sigma_{1,n}+\Sigma_{2,n}+n\Sigma_{g,n}+\dfrac{\Sigma_{h,n}}{2}\right) \notag \\
 &=-\dfrac{3h_1}{4}-h_1\left(\Sigma_{1,n}+\Sigma_{2,n}+n\Sigma_{g,n}+\dfrac{\Sigma_{h,n}}{2}\right) \label{eqNablauntn2''}.
 \end{align}

\subsection{Remainders and proof of Theorem~$\ref{thNablaun}$}

It remains to estimate $\Sigma_{1,n}, \Sigma_{2, n}, \Sigma_{g, n}$, and $\Sigma_{h,n}$
in~\eqref{eqNablauntn2''}.

\begin{lemma}\label{lSigma1}
$|\Sigma_{1,n}|\le 2EK_{h1}n^{-1/2}+\left( \dfrac{2E B_1
T_{3/2}}{a^{3/2}}+2 E^2B_2^*
\dfrac{T_2}{(\amin_n)^2}\right)n^{-1}$.
\end{lemma}
\proof Using~\eqref{eqSigma1}, we write $\Sigma_{1,n}=\Sigma'_{1,n}+\Sigma''_{1,n}, $ where
$$
\begin{aligned}
\Sigma'_{1,n}&:= \sum\limits_{|k|\le n-1}\nabla\dot y_{n-k}(a(n^2-k^2))(q_n-q_k),\\
\Sigma''_{1,n}&:= \dfrac{1}{2}\sum\limits_{|k|\le n-1}\nabla\ddot y_{n-k}(a(n^2-k^2)+\xi_{n, k})(q_n-q_k)^2.
\end{aligned}
$$

Using~\eqref{eqw01}, \eqref{eqB01}, \eqref{eqConstS11},
\eqref{eqTalpha},  and the inequality $|q_n-q_k|\le 2E n^{1/2}$,
we have
$$
\begin{aligned}
|\Sigma'_{1,n}|&\le 2En^{-1/2}\left| \sum\limits_{|k|\le n-1}\dfrac{1}{n}\dfrac{1}{a(1-x_{n, k}^2)}h'\left(\dfrac{1}{\sqrt{a}}\sqrt{\dfrac{1-x_{n, k}}{1+x_{n, k}}}\right)\right|  \\
& \;\;\;\; +\sum\limits_{|k|\le
n-1}\dfrac{2E B_1  n^{1/2}}{(a(n^2-k^2))^{3/2}} \le 2EK_{h1}n^{-1/2}+ \dfrac{2E B_1 T_{3/2}}{a^{3/2}} n^{-1}.
\end{aligned}
$$

Further, using~\eqref{eqAB2*}, \eqref{eqaminmax} and the inequalities $|\xi_{n, k}| \le |q_n-q_k|\le 2E n^{1/2}$,
we have
\begin{equation}\label{eqT1nSmall}
|\Sigma''_{1,n}|\le \dfrac{1}{2} \sum\limits_{|k|\le n-1} \dfrac{4E^2
B_2^*  n}{\left(\amin_n (n^2-k^2)\right)^2}\le  \dfrac{2 E^2 B_2^*
T_2}{(\amin_n)^2} n^{-1}.
\end{equation}
\endproof

\begin{lemma}\label{lSigma2}
$ |\Sigma_{2,n}|\le  \sum\limits_{|k|\le n-1}
\dfrac{B_0}{a(n^2-k^2)}$.
\end{lemma}
\proof Lemma \ref{lSigma2} is a straightforward consequence of~\eqref{eqSigma2} and \eqref{eqB01}.
\endproof

\begin{remark}
For $n\ge 4$, we have
$$
\sum\limits_{|k|\le n-1} \dfrac{1}{n^2-k^2}=n^{-2}+2
\sum\limits_{k=1}^{n-1}
\dfrac{1}{n^2-k^2}
\le n^{-2}+2n^{-1}\sum\limits_{k=1}^{n-1} \dfrac{1}{n-k}\le
n^{-1}(2\ln(n-1)+2+n^{-1})\to 0.
$$
\end{remark}

Propositions \ref{lemIntSqrt1} (item 1) and \ref{lemIntSqrt2} imply the following.
\begin{lemma}\label{lSigmag}
There exists $K_g > 0$ such that for $n \geq N$ the following holds\textrm{$:$} $|\Sigma_{g,n}|\le K_g n^{-3/2}$.
\end{lemma}
\begin{lemma}\label{lSigmah}
There exists $K_h > 0$ such that for $n \geq N$ the following holds\textrm{$:$} $|\Sigma_{h,n}|\le  K_h n^{-1/2}$.
\end{lemma}

Now Theorem~\ref{thNablaun} follows from~\eqref{eqNablauntn2''}
and Lemmas~\ref{lSigma1}--\ref{lSigmah}, if
the following is satisfied.
\begin{requirement}\label{ass3}
For $n \geq n_0$, the following holds\textrm{$:$}
\begin{equation}\label{eqRemainderNablaSmall}
 2EK_{h1}n^{-1/2} +\left( \dfrac{2E B_1 T_{3/2}}{a^{3/2}} + \dfrac{2 E^2 B_2^* T_2}{(\amin_n)^2} \right)n^{-1} + \sum\limits_{|k|\le n-1}
\dfrac{B_0}{a(n^2-k^2)} + \left(K_g +\dfrac{
K_h}{2}\right)n^{-1/2}\le \dfrac{3}{8}.
\end{equation}
\end{requirement}

\section{Estimates of $u_n(t), u_{n+1}(t),\dots$ for $t\in(t_{n-1},t_n)$}\label{secEstimatesUn1}

\subsection{Uniqueness of a switching moment}

As before,   we consider the sequence $t_k$ and the functions
$v_n$ given by Theorem~\ref{thqn}. In this section, we will prove
the following result.

\begin{theorem}\label{thNoSwitchBefore}
For all $n\ge n_0+1$, we have
\begin{equation}\label{equntnLess0}
v_n(t)<0,\quad
 t\in [t_{n-1},t_n) .
\end{equation}
\end{theorem}

Fix $\theta_0 >0$, satisfying the inequality ($K_h$ is given by Lemma~$\ref{lSigmah}$)
\begin{equation}
 \theta_0  K_h <\dfrac{1}{4}.\label{eqTheta0nSmall}
\end{equation}

The proof of Theorem~\ref{thNoSwitchBefore} is based on the
following proposition.
\begin{prp}\label{lNoSwitchBefore}
Let the assumptions of Theorem \str{\ref{thMainResult}} hold. Then, for all $n\ge n_0+1$,
\begin{enumerate}
\item $\dot v_n(t)\le \dfrac{3h_1  K_h}{2}(n-1)^{-1/2}$ for
all $t\in[t_{n-1},t_{n-1}+\theta_0(n-1)^{1/2}]$,

\item $\ddot v_n(t)\ge0$ for all
$t\in[t_{n-1}+\theta_0(n-1)^{1/2},t_n]$.
\end{enumerate}
\end{prp}

We first assume that Proposition~\ref{lNoSwitchBefore} is true and prove
Theorem~\ref{thNoSwitchBefore}. The proof of Proposition \ref{lNoSwitchBefore} is given in Sections \ref{sec7.2} and \ref{sec7.3} below.

\proof[Proof of Theorem~$\ref{thNoSwitchBefore}$] By \eqref{p14s} (for $n = n_0 +1$)
and Theorem~\ref{thNablaun} (for $n\ge n_0+2$), we have
$v_n(t_{n-1})\le -3h_1/8$. Therefore, Proposition~\ref{lNoSwitchBefore} (part 1) and inequality~\eqref{eqTheta0nSmall} imply that, for
$t\in [t_{n-1},t_{n-1}+\theta_0(n-1)^{1/2}]$,
\begin{equation}\label{eqTheta0nSmall'}
v_n(t)=v_n(t_{n-1})+\int_{t_{n-1}}^t \dot v_n(s)\,ds \le -\dfrac{3h_1}{8}+\theta_0(n-1)^{1/2}\dfrac{3h_1
 K_h}{2}(n-1)^{-1/2} <0.
\end{equation}

Now, for $t\in [t_{n-1}+\theta_0(n-1)^{1/2},t_n]$,
Proposition~\ref{lNoSwitchBefore} (part 2) imply that $v_n(t)$ can
vanish no more than once, and thus, by Theorem~\ref{thqn}, this
happens no earlier than at $t=t_n$.
\endproof

As a  corollary of Theorem~\ref{thNoSwitchBefore}, we obtain the
following result.

\begin{theorem}\label{thNoSwitchBeforeAllun}
For all $n\ge n_0+1$, we have
$$
v_j(t)<0,\quad
 t\in [t_{n-1},t_n),\quad j\ge n.
$$
\end{theorem}
\proof First, we show that
\begin{equation}\label{eqNoSwitchBeforeAllun1}
\nabla v_j(t_n)<0,\quad j\ge n.
\end{equation}
To do so, we estimate $\Delta v_j(t)$ for $t\in(t_{n-1},t_n]$ and
$j\ge n$. Using~\eqref{eqvInsteadOfu}, \eqref{eqy_nGreen}, and the
fact that $\dot y_n(t)\ge 0$, we have
$$
 \Delta v_j(t)= -2c-h_1\sum\limits_{|k|\le n-1} \Delta
y_{j-k}(t-t_k) =-2c-h_1\sum\limits_{|k|\le n-1} \dot
y_{j-k}(t-t_k)\le -2c.
$$
In particular, $\Delta v_j(t_n)<0$. Together with the relations
$v_n(t_n)=0$ (Theorem~\ref{thqn}) and $v_{n+1}(t_n)<0$
(Theorem~\ref{thNablaun}), this
yields~\eqref{eqNoSwitchBeforeAllun1}.

On the other hand, \eqref{eqDotynMonotoneInN} implies that $\nabla
y_{j-k}(t-t_k)$ decreases and thus
$$
\nabla v_j(t)=-c(2j+1)-h_1\sum\limits_{|k|\le n-1} \nabla
y_{j-k}(t-t_k)
$$
increases. Together with~\eqref{eqNoSwitchBeforeAllun1}, this
yields $ \nabla v_j(t)<0$ for all $t\in(t_{n-1},t_n]$ and $j\ge
n$. Since, additionally, $v_n(t)< 0$  for all $t\in[t_{n-1},t_n)$
(Theorem~\ref{thNoSwitchBefore}), the desired result follows.
\endproof

\begin{remark}
Under the assumptions of Remark~\ref{remtnRealSwitch},
Theorem~\ref{thNoSwitchBeforeAllun} implies that $t_n$ is indeed
the switching moment of the node $u_n$.
\end{remark}

\begin{remark}
Theorem~\ref{thNoSwitchBefore} implies that the
equation~in~\eqref{eq2FixedPoint} has a {\it unique} root on the
interval $[-E\sqrt{n+1},E\sqrt{n+1}]$
\end{remark}

In the rest part of this section, we will prove
Lemma~\ref{thNoSwitchBefore}

\subsection{Proof of Proposition~$\ref{lNoSwitchBefore}$\erm{:} Part \erm{1}}\label{sec7.2}

\subsubsection{Leading order terms}

We take $\theta\in[0,\theta_0(n-1)^{1/2}]$ and set
$t=t_{n-1}+\theta\in [t_{n-1},t_{n-1}+\theta_0(n-1)^{1/2}]$.

First, we represent $\dot v_n(t_{n-1}+\theta)$,
using~\eqref{eqvInsteadOfu} and the relation $t_{-(n-1)}=t_{n-1}$,
as follows:
\begin{equation}\label{eqNoSwitchBefore1'}
\begin{aligned}
\dot v_n(t_{n-1}+\theta)&=h_1-2c-h_1\sum\limits_{|k|\le n-2}\dot
y_{n-1-k}(t_{n-1}+\theta-t_k)\\
&\;\;\;\;-h_1\sum\limits_{|k|\le n-2}  \nabla \dot
y_{n-1-k}(t_{n-1}+\theta-t_k)-h_1\big(\dot y_1(\theta)+\dot
y_{2n-1}(\theta)\big).
\end{aligned}
\end{equation}

Set
$$
m:=n-1.
$$
Since $\dot y_1(\theta)\ge 0$ and $\dot
y_{2n-1}(\theta)\ge 0$ for all $\theta\ge0$, we obtain
\begin{equation}\label{eqNoSwitchBefore1}
\dot v_n(t_m+\theta) \le h_1-2c-h_1\sum\limits_{|k|\le m-1}\dot
y_{m-k}(t_m+\theta-t_k) -h_1\Sigma_{3,m},
\end{equation}
where
\begin{equation}\label{eqSigma3}
\Sigma_{3,m}:=\sum\limits_{|k|\le m-1}  \nabla \dot
y_{m-k}(t_m+\theta-t_k).
\end{equation}

Further, to apply the Taylor expansion
in~\eqref{eqNoSwitchBefore1}, we note that
$$
t_m-t_k+\theta=a(m^2-k^2)+\big(q_m-q_k+\theta\big).
$$
Therefore, using~\eqref{eqr12}, the function $H(x)=H(a,x)$ given
by~\eqref{eqHH1}, and~\eqref{eqSigmah}, we obtain
\begin{equation}\label{eqNoSwitchBefore2}
\begin{aligned}
&\sum\limits_{|k|\le m-1}\dot
y_{m-k}(t_m+\theta-t_k)=\sum\limits_{|k|\le m-1}\dot
y_{m-k}(a(m^2-k^2))+ \Sigma_{4,m}\\
&\qquad=\sum\limits_{|k|\le m-1}\dfrac{1}{m}
H(x_{m, k})+\Sigma_{5,m}+\Sigma_{4,m} =
I_H+\Sigma_{h,m}+\Sigma_{5,m}+\Sigma_{4,m},
\end{aligned}
\end{equation}
where
\begin{align}
\Sigma_{4,m}&:=\sum\limits_{|k|\le m-1}\ddot
y_{m-k}(a(m^2-k^2)+\xi_{m, k}) \big(q_m-q_k+\theta\big), \label{eqSigma4} \\
\Sigma_{5,m}&:=\sum\limits_{|k|\le m-1}
r_1(m-k,a(m^2-k^2)),\label{eqSigma5}
\end{align}
$x_{m, k}=m/k$, and $-2Em^{1/2}\le \xi_{m, k}\le 2Em^{1/2}+\theta_0
m^{1/2}$.

Combining~\eqref{eqNoSwitchBefore1}, \eqref{eqNoSwitchBefore2},
and Lemmas~\ref{lema} and~\ref{lSigmah}, we have
\begin{equation}\label{eqNoSwitchBefore3}
\dot v_n(t_m+\theta)\le h_1\left( K_h
m^{-1/2}+|\Sigma_{3,m}|+|\Sigma_{4,m}|+|\Sigma_{5,m}|\right).
\end{equation}

\subsubsection{Remainders}

It remains to estimate $\Sigma_{3,m}$, $\Sigma_{4,m}$, and $\Sigma_{5,m}$ in~\eqref{eqNoSwitchBefore3}.

\begin{lemma}\label{lSigma3}
$|\Sigma_{3,m}|\le
 K_{h1}m^{-1}+\left(  \dfrac{B_1 T_{3/2}}{a^{3/2}} +(2E+\theta_0) B_2^* \dfrac{T_2}{(\amin_m)^2}\right)m^{-3/2}$.
\end{lemma}
\proof Using~\eqref{eqSigma3}, we write $
\Sigma_{3,m}=\Sigma'_{3,m}+\Sigma''_{3,m}, $ where
$$
\begin{aligned}
\Sigma'_{3,m}&:= \sum\limits_{|k|\le m-1}\nabla\dot y_{m-k}(a(m^2-k^2)),\\
\Sigma''_{3,m}&:= \sum\limits_{|k|\le m-1}\nabla\ddot y_{m-k}(a(m^2-k^2)+\xi_{m, k})(q_m-q_k+\theta).
\end{aligned}
$$

Using~\eqref{eqw01}, \eqref{eqB01}, \eqref{eqConstS11}, and
\eqref{eqTalpha}, we have
$$
\begin{aligned}
|\Sigma'_{3,m}|&\le m^{-1}\left| \sum\limits_{|k|\le
m-1}\dfrac{1}{m}\dfrac{1}{a(1-x_{m, k}^2)}h'\left(\dfrac{1}{\sqrt{a}}\sqrt{\dfrac{1-x_{m, k}}{1+x_{m, k}}}\right)\right|\\
&+ \sum\limits_{|k|\le m-1}\dfrac{B_1}{(a(m^2-k^2))^{3/2}}\le
K_{h1}m^{-1}+ B_1  \dfrac{T_{3/2}}{a^{3/2}} m^{-3/2}.
\end{aligned}
$$

Further, using~\eqref{eqAB2*}, \eqref{eqTalpha}, and the
inequalities $|q_m-q_k+\theta|\le (2E+\theta_0) m^{1/2}$ and
$a(m^2-k^2)+\xi_{m, k}\ge \amin_m(m^2-k^2)$, we have
$$
|\Sigma''_{3,m}|\le   \sum\limits_{|k|\le m-1}
\dfrac{(2E+\theta_0)m^{1/2}B_2^*}{\left(\amin_m
(m^2-k^2)\right)^2}\le (2E+\theta_0)B_2^* \dfrac{T_2}{(\amin_m)^2}
m^{-3/2}.
$$
\endproof

\begin{lemma}\label{lSigma4}
$|\Sigma_{4,m}|\le  \dfrac{(2E+\theta_0)A_2^*
T_{3/2}}{(\amin_m)^{3/2}} m^{-1}$.
\end{lemma}
\begin{proof} Using~\eqref{eqSigma4}, \eqref{eqAB2*}, \eqref{eqTalpha},
and the inequalities $|q_m-q_k+\theta|\le (2E+\theta_0) m^{1/2}$
and $a(m^2-k^2)+\xi_{m, k}\ge \amin_m(m^2-k^2)$, we have
$$
|\Sigma_{4,m}|\le   \sum\limits_{|k|\le m-1}
\dfrac{(2E+\theta_0)A_2^* m^{1/2}}{\left(\amin_m
(m^2-k^2)\right)^{3/2}}\le
  \dfrac{(2E+\theta_0)A_2^* T_{3/2}}{(\amin_m)^{3/2}} m^{-1}.
$$
\end{proof}

\begin{lemma}\label{lSigma5}
$|\Sigma_{5,m}|\le  \dfrac{A_1 T_{3/2}}{a^{3/2}} m^{-3/2} $.
\end{lemma}
\begin{proof}
The assertion is a straightforward consequence of~\eqref{eqSigma5}, \eqref{eqA}, and~\eqref{eqTalpha}.
\end{proof}

Now inequality~\eqref{eqNoSwitchBefore3} together with
Lemmas~\ref{lSigma3}--\ref{lSigma5} yield part 1 in
Proposition~\ref{lNoSwitchBefore}, if the following requirement is satisfied
\begin{requirement}\label{ass4}
For $n \geq n_0$, the following holds\textrm{$:$}
\begin{multline}\label{eqRemainderdot1Small}
\left(K_{h1}+  \dfrac{(2E+\theta_0)A_2^*
T_{3/2}}{(\amin_m)^{3/2}}\right) m^{-1/2}
+ \left(  \dfrac{B_1
T_{3/2}}{a^{3/2}} + \dfrac{(2E+\theta_0) B_2^* T_2}{(\amin_m)^2}+
\dfrac{A_1 T_{3/2}}{a^{3/2}} \right)m^{-1}\le \dfrac{
K_h}{2}.
\end{multline}
\end{requirement}

\subsection{Proof of Proposition~$\ref{lNoSwitchBefore}$\erm{:} Part \erm{2}}\label{sec7.3}

\subsubsection{Preliminaries}

We introduce the function
\begin{equation}\label{eqPsih''h''}
\Psi(x):=h''\left(\dfrac{1}{\sqrt{a}}\sqrt{\dfrac{1-x}{1+x}}\right)+h''\left(\dfrac{1}{\sqrt{a}}\sqrt{\dfrac{1+x}{1-x}}\right),\quad
x\in(0,1).
\end{equation}
Note that
\begin{equation}\label{eqpsilim}
\lim\limits_{x\to 1}\Psi(x)=h''(0)=-\frac{1}{4\sqrt{\pi}}.
\end{equation}
Fix
\begin{equation}
\label{eqeta}
\eta \in \left(0, \frac{1}{4\sqrt{\pi}}\right).
\end{equation}

We will need the following lemma.
\begin{lemma}\label{lemsec7ep0}
There exist $x_0 \in (0, 1)$ and $\varepsilon_0>0$ such that
\begin{equation}\label{eqPsiEpsilon0Negativex0}
h''\left(\dfrac{1-x}{\sqrt{a(1-x^2)+\varepsilon_1}}\right)+h''\left(\dfrac{1+x}{\sqrt{a(1-x^2)+\varepsilon_1}}\right)\le
-\eta
\end{equation}
for all $x\in[x_0,(m-1)/m]$ and $\varepsilon_1$  with
\begin{equation}\label{eqeps0Biggerm}
-\min\left(\varepsilon_0,
 \dfrac{a}{2m}\left(2-\dfrac{1}{m}\right)\right)\le
\varepsilon_1\le\varepsilon_0.
\end{equation}
\end{lemma}

\begin{proof}
Choose $\eta_1 \in (\eta, 1/(4\sqrt{\pi}))$.
Equation \eqref{eqpsilim} implies that there is $x_0\in[0,1)$ such that
\begin{equation}\label{PsiNegativex0}
\Psi(x)\le -\eta_1,\quad x\in[x_0,1).
\end{equation}

It is not difficult to check that there exists $\varepsilon_0>0$ such
that
\begin{equation}\label{eqLimitPsiEpsilon0Negativex0}
\left|h''\left(\dfrac{1-x}{\sqrt{a(1-x^2)+\varepsilon}}\right)-
h''\left(\dfrac{1-x}{\sqrt{a(1-x^2)}}\right)\right|\le \frac{\eta_1 - \eta}{2}
\end{equation}
for all $x\in(-1,1)$ and $\varepsilon$ satisfying $|\varepsilon|\le\varepsilon_0$ and
$\varepsilon\ge -a(1-x^2)/2$.

Formulas \eqref{PsiNegativex0} and \eqref{eqLimitPsiEpsilon0Negativex0} imply Lemma~\ref{lemsec7ep0}.
\end{proof}


We fix $x_0$ and $\varepsilon_0$ from Lemma~\ref{lemsec7ep0} and  $b$ such that
\begin{equation}\label{eqbBigger2a}
b>2a.
\end{equation}

We introduce numbers $R_1$ and $R_2$ satisfying
\begin{equation}\label{eqI11}
\sum\limits_{|k|< x_0 m} \dfrac{1}{  (m^2-k^2)^{3/2}}\le
R_1m^{-2}, \quad m\ge N,
\end{equation}
\begin{equation}\label{eqI12}
 \sum\limits_{k\in[x_0m,m-1]}
\dfrac{1}{(\amax_m (m^2-k^2)+bm)^{3/2}}\ge R_2 m^{-3/2},
\quad m\ge  \max\left(N,\dfrac{1}{1-x_0}\right).
\end{equation}
Set
\begin{equation}\label{eqp30s}
B_{h2} := \sup_{x \in [0, \infty)} |h''(x)| = - h''(0), \quad B_{h4} := \sup_{x \in [0, \infty)} |h''''(x)| = h''''(0).
\end{equation}

\subsubsection{Leading order terms}

As before, we assume that $m=n-1$. We take $\theta\in[\theta_0
m^{1/2},bm]$ and set $t=t_m+\theta$. Then $t\in
[t_m+\theta_0m^{1/2},t_m+bm]$, and the latter interval contains
$[t_m+\theta_0m^{1/2},t_n]$, if the following is satisfied.
\begin{requirement}\label{ass5}
For $n \geq n_0$, the following holds \textrm{$($}see~\eqref{eqbBigger2a}\textrm{$)$}\textrm{$:$}
\begin{equation}\label{eqbmContainstn}
a(n-1)^2-E(n-1)^{1/2}+b(n-1)\ge an^2+En^{1/2}.
\end{equation}
\end{requirement}

First, we represent $\ddot v_n(t_{n-1}+\theta)$,
using~\eqref{eqvInsteadOfu}, as follows
(cf.~\eqref{eqNoSwitchBefore1'}):
\begin{equation}\label{eqvnddot1}
\ddot v_n(t_m+\theta)=-h_1(I_{1,m}+I_{2,m}+\Sigma_{6,m}),
\end{equation}
where
\begin{align}
I_{1,m}&:=\sum\limits_{|k|\le m-1}\ddot
y_{m-k}(t_m+\theta-t_k),\label{eqI12ddotv} \quad
I_{2,m}:=\ddot y_1(\theta)+\ddot y_{2m+1}(\theta),\\
\Sigma_{6,m}&:=\sum\limits_{|k|\le m-1}  \nabla \ddot
y_{m-k}(t_m+\theta-t_k).\label{eqSigma6}
\end{align}

 Using~\eqref{eqr12}, we represent $I_{1,m}$ as
follows:
\begin{multline}\label{eqI1ddoty}
I_{1,m} = \sum\limits_{|k|\le m-1}
\frac{1}{(a(m^2-k^2)+q_m-q_k+\theta)^{3/2}} h''\left(
\frac{m-k}{\sqrt{a(m^2-k^2)+q_m-q_k+\theta}}\right)+\Sigma_{7,m},
\end{multline}
where
\begin{equation}\label{eqSigma7}
\Sigma_{7,m} :=\sum\limits_{|k|\le m-1}
r_2(m-k,a(m^2-k^2)+q_m-q_k+\theta).
\end{equation}

Now we split the sum in~\eqref{eqI1ddoty} into two sums in which the
summation is taken over $|k|< x_0m$ and $|k|\in[x_0m,m-1]$,
respectively. Let us estimate the first sum,
using~\eqref{eqI11}, \eqref{eqp30s}, and the inequality
$a(m^2-k^2)+q_m-q_k\ge \amin_m(m^2-k^2)$:
\begin{equation}\label{eqI1ddoty2}
\begin{aligned}
&\left|\sum\limits_{|k|< x_0 m}
\dfrac{1}{(a(m^2-k^2)+q_m-q_k+\theta)^{3/2}} h''\left(
\dfrac{m-k}{\sqrt{a(m^2-k^2)+q_m-q_k+\theta}}\right)\right|\\
&\qquad \le \sum\limits_{|k|< x_0 m} \dfrac{B_{h2}}{(\amin_m (
(m^2-k^2))^{3/2}}\le  \dfrac{B_{h2} R_1}{(\amin_m)^{3/2}}
m^{-2}.
\end{aligned}
\end{equation}

To estimate the second sum (which we do if $m-1\ge x_0 m$,
i.e., $m\ge 1/(1-x_0)$), we set
$$
\varepsilon_{m, k}:=\dfrac{q_m-q_k+\theta}{m^2}.
$$

Below we assume that the following holds.
\begin{requirement}\label{ass6.1}
For $n\ge
\max\left(n_0(E),\frac{1}{1-x_0}\right)$, the following holds\textrm{$:$}
$$
-\min\left(\varepsilon_0,
 \dfrac{a}{2m}\left(2-\dfrac{1}{m}\right)\right)\le \dfrac{(\theta_0-2E)m^{1/2}}{m^2}.
$$
\end{requirement}

\begin{requirement}\label{ass6.2}
For $n\ge
\max\left(n_0(E),\frac{1}{1-x_0}\right)$, the following holds\textrm{$:$}
$$
\dfrac{2Em^{1/2}+bm}{m^2}\le \varepsilon_0.
$$
\end{requirement}

Then
\begin{equation}\label{eqbm-1Small}
 -\min\left(\varepsilon_0, \frac{a}{2m}\left(2-\frac{1}{m}\right)\right)
 \le \varepsilon_{m, k}
 \le \varepsilon_0, \quad \quad \amin_m \le \frac{a(m^2-k^2)+q_m-q_k}{m^2-k^2}\le
\amax_m
\end{equation}
Using~\eqref{eqbm-1Small} as well
as~\eqref{eqPsiEpsilon0Negativex0} and~\eqref{eqI12}, we have ($x_{m, k}$ is given by \eqref{eqxnk})
\begin{equation}\label{eqI1ddoty3}
\begin{aligned}
&\sum\limits_{|k|\in[x_0m,m-1]}\!
\frac{1}{(a(m^2-k^2)+q_m-q_k+\theta)^{3/2}} h''\!\!\left(\!
\frac{m-k}{\sqrt{a(m^2-k^2)+q_m-q_k+\theta}}\!\right)\\
&=\sum\limits_{k\in[x_0m,m-1]}
\frac{1}{(a(m^2-k^2)+q_m-q_k+\theta)^{3/2}}\\
& \quad \quad \cdot\left[h''\left(
 \frac{1-x_{m, k}}{\sqrt{a(1-x_{m, k}^2)+\varepsilon_{m, k}}}\right)+h''\left(
\frac{1+x_{m, k}}{\sqrt{a(1-x_{m, k}^2)+\varepsilon_{m, k}}}\right)\right]\\
&\le -\eta \sum\limits_{k\in[x_0m,m-1]}
\frac{1}{(\amax_m (m^2-k^2)+bm)^{3/2}}\le
-\eta R_2 m^{-3/2}.
\end{aligned}
\end{equation}

Thus, using~\eqref{eqI1ddoty}, \eqref{eqI1ddoty2},
and~\eqref{eqI1ddoty3}, we obtain the following estimate for
$I_{1,m}$ in~\eqref{eqI12ddotv}:
\begin{equation}\label{eqI1ddoty4}
I_{1,m} \le - \eta R_2 m^{-3/2}+
\dfrac{B_{h2} R_1}{(\amin_m)^{3/2}} m^{-2}+\Sigma_{7,m}.
\end{equation}

In what follows, we assume that the following is satisfied.
\begin{requirement}\label{asstau02}
For $n \geq n_0$, the following holds\textrm{$:$}
$$
\tau_0 \leq \theta_0 n^{1/2}.
$$
\end{requirement}

Now we represent $I_{2,m}$ in~\eqref{eqI12ddotv},
using~\eqref{eqr12} and the equalities $h''(0)=-1/(4\sqrt{\pi})$
and $h'''(0)=0$, as follows:
$$
\begin{aligned}
&I_{2,m}=\dfrac{1}{\theta^{3/2}}h''\left(\dfrac{1}{\sqrt{\theta}}\right)+
\dfrac{1}{\theta^{3/2}}h''\left(\dfrac{2m+1}{\sqrt{\theta}}\right)+ r_2(1, \theta) + r_2(2m+1, \theta)\\
&=\dfrac{1}{\theta^{3/2}}\!\left(\!-\dfrac{1}{4\sqrt{\pi}}+\dfrac{h''''(\xi)}{2\theta}+
h''\left(\dfrac{2m+1}{\sqrt{\theta}}\right)+\theta^{3/2}\bigl(
r_2(1, \theta) +  r_2(2m+1, \theta) \bigr)\!\right)\!,
\end{aligned}
$$
where $\xi\in[0,\theta^{-1/2}]$.
Below we assume that the following holds.
\begin{requirement}\label{ass7.1}
For $n \geq n_0$ and $m = n-1$, the following holds\textrm{$:$}
$$
\dfrac{B_{h4}}{2\theta_0}m^{-1/2}+
\sup_{x \geq \frac{2m+1}{\sqrt{bm}}}
h''\left(x\right)
+
\dfrac{2A_2}{\theta_0^{1/2}}m^{-1/2}\le \dfrac{1}{4\sqrt{\pi}}.
$$
\end{requirement}
Hence,
\begin{equation}\label{eqI2ddoty}
I_{2,m}\le 0.
\end{equation}

Due to~\eqref{eqvnddot1}, \eqref{eqI1ddoty4}, and
\eqref{eqI2ddoty},
Proposition~$\ref{lNoSwitchBefore}$ (part 2) follows if
\begin{equation}\label{eqvnddot2}
-\eta R_2 m^{-3/2}+ \dfrac{B_{h2}
R_1}{(\amin_m)^{3/2}}
m^{-2}+|\Sigma_{7,m}|+|\Sigma_{6,m}|\le 0.
\end{equation}

\subsubsection{Remainders}
Let us prove~\eqref{eqvnddot2}. To do so, we need to estimate
$\Sigma_{6,m}$ and $\Sigma_{7,m}$.
\begin{lemma}\label{lSigma6}
$|\Sigma_{6,m}| \le \dfrac{B_2^* T_2}{(\amin_m)^2} m^{-2}  $.
\end{lemma}
\proof Using~\eqref{eqSigma6}, \eqref{eqAB2*}, \eqref{eqTalpha}, and \eqref{eqbm-1Small}, we have
$$
|\Sigma_{6,m}|\le \sum\limits_{|k|\le m-1} \dfrac{B_2^*}{(\amin_m
(m^2-k^2))^2}\le \dfrac{B_2^* T_2}{(\amin_m)^2} m^{-2}.
$$

\begin{lemma}\label{lSigma7}
$|\Sigma_{7,m}|\le \dfrac{A_2T_{5/2}}{(\amin_m)^{5/2}} m^{-2}  $.
\end{lemma}
\proof Using~\eqref{eqSigma7}, \eqref{eqA}, \eqref{eqTalpha}, and \eqref{eqbm-1Small}, we have
$$
|\Sigma_{7,m}|\le \sum\limits_{|k|\le m-1} \dfrac{A_2}{(\amin_m
(m^2-k^2)+\theta_0 m^{1/2})^{5/2}}\le
\dfrac{A_2T_{5/2}}{(\amin_m)^{5/2}} m^{-5/2}.
$$

Using Lemmas~\ref{lSigma6} and~\ref{lSigma7}, we see
that~\eqref{eqvnddot2} holds, if the following is satisfied.
\begin{requirement}\label{ass8}
For $n \geq n_0$ and $m = n-1$, the following holds\textrm{$:$}
\begin{equation}\label{eqvnddot2'}
\left( \dfrac{B_{h2} R_1}{(\amin_m)^{3/2}}+ \dfrac{B_2^*
T_2}{(\amin_m)^2}\right)m^{-1/2}+\dfrac{A_2T_{5/2}}{(\amin_m)^{5/2}}m^{-1}\le
\eta R_2.
\end{equation}
\end{requirement}

\section{Main result: proof of Theorem~\ref{thMainResult}}\label{secProofMainTheorem}

For $n=n_0+1$, Theorems~\ref{thqn} and \ref{thNoSwitchBeforeAllun} imply that the node $u_n(t)$
achieves the threshold $0$ at a time moment $t_n=an^2+q_n$, where
$q_n\in[-E\sqrt{n},0]$ with the same $E$ as in \eqref{eqEn+ScondtkForuntn}.
Moreover, neither of the nodes $u_n,u_{n+1},\dots$ switches on the interval
$[t_{n-1},t_n)$, and thus $u_n(t)$ switches exactly at the moment
$t_n$. Furthermore, by Theorem~\ref{thNablaun}, the required
estimates for $\nabla u_n(t_n)$ hold. Thus, the assertion of
Theorem~\ref{thMainResult} holds for $n=n_0+1$.

In particular, we see that items 1 and 2 in Definition~\ref{condtkForuntn} hold
with $n_0$ replaced by $n_0+1$ and with the same $E$ as before.
Hence, we can repeat the above argument to obtain the assertion of
Theorem~\ref{thMainResult} for $n=n_0+2$,  and so on, by
induction, for any $n\ge n_0+1$.

\appendix

\section{Equivalence of three equations: proof of Proposition~\ref{lema}}
\label{secaInt}

We prove that
\begin{enumerate}
\item equation \eqref{eqah1h} has a unique root,
\item equations \eqref{eqah1f} and \eqref{eqah1h} have the same roots,
\item equations \eqref{eqah1g} and \eqref{eqah1h} have the same roots.
\end{enumerate}

We will prove in detail items 1 and 2.


Making the change of variables $y = \dfrac{1}{\sqrt{a}}\sqrt{\dfrac{1-x}{1+x}}$ in \eqref{eqHH1} and \eqref{eqFG}, we have
\begin{equation}\notag
I_H(a)  = \int_0^{+\infty} \frac{2h(y)}{1+ ay^2} dy, \quad I_G(a) = \int_0^{+\infty} \frac{4ayg(y)}{(1+ ay^2)^2} dy, \quad I_F(a)  = \int_0^{+\infty} \frac{8a^2y^2f(y)}{(1+ ay^2)^3} dy.
\end{equation}
Now we see that $I_H(a)$ decreases from 1 to 0 as $a$ increases from $0$ to $+\infty$. Hence, for any $0 < c <h_1/2$, equation \eqref{eqah1h} has a unique root $a > 0$. Item 1 is proved.

Let us prove item 2.
Integrating by parts and using equations \eqref{eqfghft} and \eqref{eqfgh}, we obtain
\begingroup
\begin{align*}
I_F(a) & = \int_0^{+\infty} \left( \frac{8a^2y^2}{(1+ ay^2)^3} - \frac{a}{1+ay^2} \right) f(y) dy + \int_0^{+\infty} \frac{a}{1+ay^2} f(y) dy \\
& =  \left.\frac{ay(-1 + ay^2)}{(1+ay^2)^2}g(y)\right|_0^{+\infty} - \int_0^{+\infty} \frac{ay(-1 + ay^2)}{(1+ay^2)^2}g(y) dy + \int_0^{+\infty} \frac{a}{1+ay^2} f(y) dy  \\
& = \int_0^{+\infty} \left( - \frac{ay(-1 + ay^2)}{(1+ay^2)^2} + y \frac{a}{1+ay^2} \right) g(y) dy + \int_0^{+\infty} \frac{a}{1+ay^2}\left(f(y) - y g(y)\right) dy \\
& = - \left.\frac{1}{1+ ay^2}g(y)\right|_0^{+\infty} + \int_0^{+\infty} \frac{1}{1+ ay^2}h(y) dy + \int_0^{+\infty} \frac{2a}{1+ ay^2}h(y) dy \\
& = - \frac{1}{2} + \int_0^{+\infty} \frac{2a+1}{1+ay^2} h(y) dy = \frac{2a+1}{2}I_H(a) - \frac{1}{2}.
\end{align*}
\endgroup

It is easy to conclude from the last identity that equations \eqref{eqah1f} and \eqref{eqah1h} have the same roots. Item 2 is proved.

Similarly to item 2, integrating by parts and using relations \eqref{eqfghft} and \eqref{eqfgh}, we obtain
$$
I_G(a) = \int_0^{+\infty} \frac{4ay}{(1+ ay^2)^2}g(y) dy = \left. -\frac{2}{1+ay^2}g(y)\right|_0^{+\infty} + \int_0^{+\infty}\frac{2}{1+ay^2}h(y) = -1 + I_H(a).
$$
From this identity, it is easy to conclude item 3.



\section{Requirements on $n_0(E)$}

In this appendix, we collect the constants that we use throughout
the paper to determine $n_0(E)$ as well as all the 12 requirements on
the number $n_0=n_0(E)$ entering Definition~\ref{condtkForuntn} of admissible $E$.

\subsection{Constants not depending on $a$ or $E$} \label{subsecConstNoAE}

\begin{enumerate}
\item $\tau_0 > 0$ is an arbitrarily fixed real number (see \eqref{eqtau0}).

\item Set (see \eqref{eqA})
$
A_0 := \sup_{n \geq 0, t \geq \tau_0} \sqrt{t}\lra{y_n(t)  -
\sqrt{t}\,f\left(\frac{n}{\sqrt{t}}\right)}.$

\item Set (see \eqref{eqA})
$
A_1 :=
\sup_{n \geq 0, t \geq \tau_0} t\sqrt{t}\lra{\dot y_n(t)  -
\frac{1}{\sqrt{t}}\,h\left(\frac{n}{\sqrt{t}}\right)}.
$

\item Set (see \eqref{eqA})
$
\At_1 := \sup_{n \geq 0, t \geq \tau_0} t\sqrt{t}\lra{y_n(t)  - \sqrt{t}\,f\left(\frac{n}{\sqrt{t}}\right) - \frac{1}{\sqrt{t}}\,\ft\left(\frac{n}{\sqrt{t}}\right)}.
$

\item\label{itemBX} Set (see \eqref{eqA})
$
A_2 :=
\sup_{n \geq 0, t \geq \tau_0} t^2\sqrt{t} \lra{\ddot y_n(t)  -
\frac{1}{t\sqrt{t}}\,h''\left(\frac{n}{\sqrt{t}}\right)}.
$

\item Set (see \eqref{eqB01})
$
B_0 := \sup_{n \geq 0, t \geq \tau_0} t\lra{\nabla y_n(t)  - g\left(\frac{n}{\sqrt{t}}\right) -
\frac{1}{2\sqrt{t}}\,h\left(\frac{n}{\sqrt{t}}\right)}.
$

\item Set (see \eqref{eqB01})
$
B_1 := \sup_{n \geq 0, t \geq \tau_0} t\sqrt{t}\lra{
\nabla \dot y_n(t) - \frac{1}{t}h'\left(\frac{n}{\sqrt{t}}\right)}.
$

\item Set (see \eqref{eqAB2*})
$
A_2^* :=
\sup_{n \geq 0, t \geq \tau_0} t\sqrt{t} |\ddot y_n(t)|.
$

\item Set (see \eqref{eqAB2*})
$
B_2^* := \sup_{n \geq 0, t \geq \tau_0} t^2 |\nabla\ddot y_n(t)|.
$

\item We use the notation (see \eqref{eqR})
$
R := \pi.
$

\item $B_{h2}$ is given by (see \eqref{eqp30s})
$
B_{h2}:=\sup\limits_{x\in[0,\infty)}
\left|h''(x)\right|=-h''(0).
$

\item $B_{h4}$ is given by (see \eqref{eqp30s})
$
B_{h4}:=\sup\limits_{x\in[0,\infty)}
\left|h''''(x)\right|=h''''(0).
$


\end{enumerate}

\subsection{Constants depending on $a$ but not depending on $E$}

\begin{enumerate}
\item $a$ is a unique root of equation~\eqref{eqah1f} (or equivalently \eqref{eqah1g}, \eqref{eqah1h}).

\item Consider the function
$
h_a(x) := h\lr{\frac{x}{\sqrt{a}}} + h\lr{\frac{1}{x\sqrt{a}}}.
$
Set (see \eqref{eqfDa})
$
D_a := \inf_{x \in (0, 1]} h_a(x).
$

\item Set (see \eqref{eqp})
$
p  := \sup_{x \in (0, 1]} \frac{h_a'(x)x}{h_a(x)}.
$

\item $N$ is a fixed natural number satisfying (see \eqref{eqep}, \eqref{eqkap0})\\
$
N\lr{\lr{\frac{N+1}{N}}^{\frac{1+p}{2}}-1} - \lr{2^{\frac{1-p}{2}}-1}
- \lr{2^{\frac{1-p}{2}}-1+ N\lr{\lr{\frac{N+1}{N}}^{\frac{1+p}{2}}-1}}\lr{2^{\frac{1-p}{2}}-1}\frac{2}{N} > 0.
$

\item Set (see \eqref{eqDp12})
$
D_{p1}  := 2^{\frac{1-p}{2}}-1$,
$D_{p2} := N\lr{\lr{1+\frac{1}{N}}^{\frac{1+p}{2}}-1}$.

\item Constants needed to define $K$ and $K'$ from Lemma \ref{lemCn}
    \begin{enumerate}
    \item Set (see \sref{eqK1p})
    $
    L_1(F, N) := \sup_{n \geq N} n^{3/2} \left|\int_{-1}^1 F(x) dx -  \sum_{ |k| \leq n-1} \frac{1}{n}F\lr{\frac{k}{n}} \right|,
    $
    where the supremum exists due to Proposition \ref{lemIntSqrt1}. Set
    $
    K'_1  := h_1 L_1(F, N).
    $
    \item Set (see \eqref{eql1})
    $
    l_1 := h_1A_0\frac{1}{a^{1/2}}R.
    $
    \item Set (see \eqref{eqK'})
    $
    K' := K'_1 + l_1 \frac{1}{\sqrt{N}}.
    $
    \item Set
    $
    K_1 := \sup_{n \geq N} n^{1/2}|C_{n+1}^{(1)} - C_n^{(1)}|,
    $
    where
     $$
     C_l^{(1)}  := h_1\lr{l^2 I_F - \sum_{|k| \leq l-1} \sqrt{a(l^2-k^2)}f\lr{\frac{l-k}{\sqrt{a(l^2-k^2)}}}}, \quad l = n, n+1.
     $$
     Note that the supremum exists due to \eqref{P16s}.
    \item Consider the function
    $
    \tilde{F}(x) := \frac{1}{\sqrt{a(1-x^2)}}\tilde{f}\lr{\frac{1-x}{\sqrt{a(1-x^2)}}}
    $
    and constants $L_{\ft}$, $L_{\ft}^*, l_{\ft} \geq 0$ such that (see \eqref{eqLt})
    $$
    L_{\tilde{f}}^* \frac{1}{\sqrt{n}} - l_{\tilde{f}}\frac{1}{n}
    \leq
    -\int_{-1}^1 \Ft(x) \dx + \sum_{-n \leq k \leq n-1} \frac{1}{n}\Ft\lr{\frac{k}{n}}
    \leq
    L_{\tilde{f}} \frac{1}{\sqrt{n}}, \quad n \geq N.
    $$
    Set (see \eqref{P16ss})
    $
    K_2 := h_1 \lr{ L_{\tilde{f}} - L_{\tilde{f}}^* + l_{\tilde{f}}\frac{1}{\sqrt{N}}}.
    $
    \item Set (see \eqref{P16sss})
    $
    K_3 := 2h_1\tilde{A}_1 \frac{1}{a^{3/2}}T_{3/2}.
    $
    \item Set (see \eqref{eqK})
    $
    K := K_1 + K_2 + K_3 \frac{1}{N}.
    $
    \end{enumerate}

\begin{remark}
In principle, due to Proposition \ref{lemCn}, we could define $K' := \sup_{n \geq N} n^{-1/2}|C_n|$ and $K := \sup_{n \geq N} n^{1/2}|C_{n+1}-C_n|$. However, calculation of the values $C_n$ is computationally consuming as it involves Bessel functions. We used the strategy described above, since it is based  on error estimates of Riemann sums only.
\end{remark}

\item Set (see \eqref{eqep})
$
\varkappa  := D_{p2} - D_{p1} - 2(D_{p1}+ D_{p2})D_{p1}\frac{1}{N} \; (>0).\\
$
\item Set (see \eqref{eqE0})
$
E_0 := \frac{K + \varkappa K'}{(h_1 - 2c)\varkappa} \; (>0).
$
\item For $\al> 0$, set (see \sref{eqSal})
$
S_{\al} := \sup_{n \geq N} \lr{\sqrt{n}\sum_{|k| \leq n-1}\frac{(n-k)^{\al}}{(n^2-k^2)^{\al+1/2}}}.
$
We use only the values of $S_1$, $S_2$, and $S_3$ to determine $n_0(E)$.
\item For $\al > 1$, set (see \sref{eqTalpha})
$
T_\alpha := \sup_{n\geq N}\lr{  \sum_{|k|\le n-1}\frac{n^{\al}}{(n^2-k^2)^{\al}}}.
$
We use only the values of $T_{3/2}$, $T_2$, and $T_{5/2}$ to determine
$n_0(E)$.
\item Let constants $L_2^*(H, N) > \frac{1}{\sqrt{2a}}h(0)$ and $l_2(H, N) \geq 0$ be such that (see Proposition~\ref{lemIntSqrt2})
$$
L_2^*(H, N)\frac{1}{n^{1/2}} - l_2(H, N) \frac{1}{n} \leq I_H -
\sum_{ k = -n}^{n-1} \frac{1}{n}H\lr{\frac{k}{n}}, \quad n \geq N,
$$
where $H(x)$ is given by \eqref{eqHH1} and $I_H = I_H(a)$ is given by \eqref{eqIntfgh}.
Set (see \eqref{eqCH})
$C_H  := L_2^*(H, N)$, $l_H  := l_2(H, N)$.
\item Consider the function $\bar{H}(x) := \frac{1}{1-x^2}H(x)$. Set (see \eqref{eqCH2})\\
$
C_{H2} := \sup_{n \geq N} n^{-1/2} \left| \sum_{|k| \leq n-1}
\frac{1}{n} \bar{H}\lr{\frac{k}{n}} \right|.
$
\item Set (see~\eqref{eqConstS11})
$
K_{h1} := \sup_{n \geq N} n^{-1}\left| \sum\limits_{|k|\le
n-1}\frac{1}{a(1-(k/n)^2)}h'\left(\frac{1}{\sqrt{a}}\sqrt{\frac{1-k/n}{1+k/n}}\right)\right|.
$
\item Set (see Lemma~\ref{lSigmag} and \eqref{eqSigmag})
$
K_g := \sup_{n \geq N} n^{3/2}\left|\sum\limits_{|k|\le
n-1}\frac{1}{n} G\lr{\frac{k}{n}}-I_G-\frac{1}{4n}\right|,
$
where
$G(x)$ is given by~\eqref{eqFG}, and $I_G=I_G(a)$  is given
by~\eqref{eqIntfgh}.
\item Set (see Lemma~\ref{lSigmah} and \eqref{eqSigmah})
$
K_h := \sup_{n \geq N} n^{1/2}\left|\sum\limits_{|k|\le
n-1}\frac{1}{n} H\lr{\frac{k}{n}}-I_H \right|,
$
where
$H(x)$ is given by~\eqref{eqHH1}, and $I_H=I_H(a)$ is given by~\eqref{eqIntfgh}.
\item $\theta_0$ satisfies (see~\eqref{eqTheta0nSmall})
$
 \theta_0
 K_h <\frac{1}{4}.
$
\item $\eta$ satisfies (see \eqref{eqeta})
$
\eta \in \left(0, \frac{1}{4\sqrt{\pi}}\right).
$
\item $x_0\in[0,1)$ and $\varepsilon_0>0$ satisfy
(see~\eqref{eqPsiEpsilon0Negativex0}, \eqref{eqeps0Biggerm})
$$
h''\left(\frac{1-x}{\sqrt{a(1-x^2)+\varepsilon_1}}\right)+h''\left(\frac{1+x}{\sqrt{a(1-x^2)+\varepsilon_1}}\right)\le
-\eta
$$
for all $x\in[x_0,(n-1)/n]$ and $\varepsilon_1$  with
$$
-\min\left(\varepsilon_0,
 \frac{a}{2n}\left(2-\frac{1}{n}\right)\right)\le
\varepsilon_1\le\varepsilon_0,\quad n\ge
\max\left(N,\frac{1}{1-x_0}\right)
$$
\item $b$ satisfies (see~\eqref{eqbBigger2a})
$
b>2a.
$
\item Set (see~\eqref{eqI11})
$
R_1 := \sup_{n \geq N} n^2\sum\limits_{|k|< x_0 n} \frac{1}{ (n^2-k^2)^{3/2}}
$
\end{enumerate}

\subsection{Constants depending on $E$}\label{subsecConstE}

\begin{enumerate}
\item
Set (see \eqref{eqaminmax})
$
\amin_n:=a-\frac{2E n^{1/2}}{2n-1}.
$
\item
Set (see \eqref{eqaminmax})
$
\amax_n:=a+\frac{2E n^{1/2}}{2n-1}.
$
\item Set (see \eqref{eqdeltan})
$
\delta_n := \frac{2 A_1}{aD_a}\frac{1}{n} + \frac{ 2A_2^*E
a^{1/2}}{D_a (\amin_n)^{3/2}}
\frac{1}{\sqrt{n}}.
$
\item Set  (see~\eqref{eqI12})
$
R_2 := \inf_{n\ge \max\left(N,1/(1-x_0)\right)}  \sum\limits_{k\in[x_0n,n-1]}
\frac{n^{3/2}}{(\amax_n(n^2-k^2)+bn)^{3/2}}.
$
\end{enumerate}

\subsection{Requirements on $n_0(E)$}\label{refAppendixSubsectionAssumptionsN0}

We assume that the following requirements hold for $n\ge n_0(E)$:

\begin{enumerate}
\item
$n\ge N.$
\item
$
\amin_n \geq \frac{\tau_0}{2n-1}.
$
\item
$
\delta_n + \delta_{n+1} \leq 2D_{p1}.
$
\item
$
C_H - \frac{h(0)}{\sqrt{2a}} \geq
\frac{l_H}{\sqrt{n}} + \frac{\varkappa h(0)}{\sqrt{2a}\cdot n}  + C_{H2}\delta_n +  \frac{1}{\sqrt{2a}}\!\left(\!1 \! + \! \frac{\varkappa}{n}\!\right)\! 
\left(\!
h(0)\frac{n+1}{2n+1}\delta_{n+1} +
(1+\delta_{n+1})h\!\left(\!\!\sqrt{\frac{2n+1}{a}}\right)
\!\right).
$
\item
$
 2EK_{h1}n^{-\frac{1}{2}} +\left( \frac{2E B_1 T_{3/2}}{a^{3/2}} +
 \frac{2 E^2 B_2^* T_2}{(\amin_n)^2}\right)n^{-1}  +
\sum\limits_{|k|\le n-1} \frac{B_0}{a(n^2-k^2)} + \left(K_g
+\frac{ K_h}{2}\right)n^{-\frac{1}{2}}\le \frac{3}{8}.
$
\item
$
\left(K_{h1}+
 \frac{(2E+\theta_0)A_2^* T_{3/2}}{(\amin_n)^{3/2}} \right)
n^{-1/2}+\left(  \frac{B_1 T_{3/2}}{a^{3/2}} +
\frac{(2E+\theta_0) B_2^* T_2}{(\amin_n)^2}+
 \frac{A_1 T_{3/2}}{a^{3/2}} \right) n^{-1}\le
\frac{ K_h}{2}.
$
\item
$
a(n-1)^2-E(n-1)^{1/2}+b(n-1)\ge an^2+En^{1/2}.
$
\item
$
\frac{2E-\theta_0 }{n^{3/2}}\le \min\left(\varepsilon_0,
 \frac{a}{2n}\left(2-\frac{1}{n}\right)\right)
$
\; for \;
$
n \ge \max\left(n_0(E),\frac{1}{1-x_0}\right).
$
\item
$
2En^{-3/2}+b n^{-1}\le \varepsilon_0
$
\; for \;
$
n\ge \max\left(n_0(E),\frac{1}{1-x_0}\right).
$
\item
$
\tau_0 \leq \theta_0 n^{1/2}.
$
\item
$
  \frac{B_{h4}}{2\theta_0}n^{-1/2}
+ \sup_{x \geq \frac{2n+1}{\sqrt{bn}}} h''\left(x\right)
+ \frac{2A_2}{\theta_0^{1/2}}n^{-1/2}\le
\frac{1}{4\sqrt{\pi}}.
$
\item
$
\left(\frac{B_{h2}  R_1}{(\amin_n)^{3/2}}+ \frac{B_2^*
T_2}{(\amin_n)^2}\right)n^{-1/2}+
\frac{A_2T_{5/2}}{(\amin_n)^{5/2}} n^{-1}\le \eta
R_2.
$
\end{enumerate}

\section*{Acknowledgement}
The authors are grateful to Daria Neverova for her help in preparing the figures.
The work of the first author was supported by the DFG
Heisenberg Programme, DFG project SFB 910, and the Ministry of Education and
Science of Russian Federation (agreement 02.a03.21.0008).
The second author would like to thank JSC ``Gazprom neft'' and Contest ``Young Russian Mathematics'' for their attention to this work.

\end{document}